\documentclass{amsart}
\usepackage{pictex}

\begin{document}

\newtheorem{thm}{Theorem}[section]
\newtheorem{lem}[thm]{Lemma}
\newtheorem{cor}[thm]{Corollary}
\newtheorem{add}[thm]{Addendum}
\newtheorem{prop}[thm]{Proposition}

\theoremstyle{definition}
\newtheorem{defn}[thm]{Definition}

\theoremstyle{remark}
\newtheorem{rmk}[thm]{Remark}

\def\square{\hfill${\vcenter{\vbox{\hrule height.4pt \hbox{\vrule width.4pt
height7pt \kern7pt \vrule width.4pt} \hrule height.4pt}}}$}

\newenvironment{pf}{{\it Proof:}\quad}{\square \vskip 12pt}

\title{Generalizations of McShane's Identity to Hyperbolic Cone-Surfaces}
\author{Ser Peow Tan, Yan Loi Wong, and Ying Zhang}
\address{Department of Mathematics \\ National University of Singapore \\ 2 Science Drive 2 \\ Singapore 117543}
\email{mattansp@nus.edu.sg; matwyl@nus.edu.sg;
scip1101@nus.edu.sg}

\maketitle
%
%

\begin{abstract}
We generalize McShane's identity for the length series of simple
closed geodesics on a cusped hyperbolic surface
\cite{mcshane1998im} to hyperbolic cone-surfaces (with all cone
angles $\le \pi$), possibly with cusps and/or geodesic boundary.
In particular, by applying the generalized identity to the
orbifolds obtained from taking the quotient of the one-holed torus
by its elliptic involution, and the closed genus two surface by
its hyper-elliptic involution, we obtain generalizations of the
Weierstrass identities for the one-holed torus, and identities for
the genus two surface, also obtained by McShane using different
methods in \cite{mcshane2004blms}, \cite{mcshane1998preprint} and
\cite{mcshane2004preprint}. We also give an interpretation of the
identity in terms of complex lengths, gaps, and the direct visual
measure of the boundary.
\end{abstract}
\section{{\bf Introduction}}\label{s:intro}
\vskip 20pt

Greg McShane discovered the following striking identity in his
Ph.D. thesis:

\begin{thm}\label{thm:mcshane torus} {\rm (McShane
\cite{mcshane1991thesis})} In a once punctured hyperbolic torus
$T$,
\begin{eqnarray}
\sum_{\gamma}\frac{1}{1 + \exp |\gamma|} = \frac{1}{2},
\end{eqnarray}
where the sum extends over all simple closed geodesics on $T$ and
where $|\gamma|$ denotes the length of $\gamma$ in the given
hyperbolic structure.
\end{thm}

Throughout this paper we shall always use $|\gamma|$ to denote the
hyperbolic length of $\gamma$ if $\gamma$ is a (generalized)
simple closed geodesic or a simple geodesic arc on a hyperbolic
(cone-)surface. All surfaces considered in this paper are assumed
to be connected and orientable.

\vskip 10pt

Later McShane extended his identity to more general surfaces:

\begin{thm}\label{thm:mcshane general} {\rm (McShane
\cite{mcshane1998im})} In a finite area hyperbolic surface $M$
with cusps and without boundary,
\begin{eqnarray}
\sum \frac{1}{1 + \exp \frac{1}{2}(|\alpha|+|\beta|)} =
\frac{1}{2},
\end{eqnarray}
where the sum is over all unordered pairs of simple closed
geodesics $\alpha, \beta$ (where $\alpha$ or $\beta$ might be a
cusp treated as a simple closed geodesic of length $0$) on $M$
such that $\alpha, \beta$ bound with a distinguished cusp point an
embedded pair of pants on $M$.
\end{thm}

Note that Theorem \ref{thm:mcshane torus} can be regarded as a
special case of Theorem \ref{thm:mcshane general} where $\alpha,
\beta$ are the same for each pair $\alpha, \beta$.

\vskip 10pt

In \cite{mcshane2004blms} McShane demonstrated three other closely
related identities for the lengths of simple closed geodesics in
each of the three Weierstrass classes on a hyperbolic torus.
Recall that a hyperbolic torus $T$ has three Weierstrass points
which are the fixed points of the unique elliptic involution which
maps each simple closed geodesic on $T$ onto itself with
orientation reversed, and for a Weierstrass point $x$ on $T$ the
simple closed geodesics in the Weierstrass class which is dual to
$x$ are precisely all the simple closed geodesics on $T$ which do
not pass through $x$.

\vskip 10pt

\begin{thm}\label{thm:mcshane weierstrass} {\rm (McShane
\cite{mcshane2004blms})} In a once punctured hyperbolic torus,
\begin{eqnarray}\label{eqn:mcshane weierstrass}
\sum_{\gamma \in \mathcal A} \sin^{-1} \left ( \frac{1}{\cosh
\frac{1}{2}|\gamma|} \right )= \frac{\pi}{2},
\end{eqnarray}
where the sum is over all simple closed geodesics in a Weierstrass
class $\mathcal A$.
\end{thm}

\vskip 10pt

On the other hand, B. H. Bowditch gave an alternative proof of
Theorem \ref{thm:mcshane torus} using Markoff triples
\cite{bowditch1996blms} and extended the identity in Theorem
\ref{thm:mcshane torus} to the case of quasi-fuchsian
representations of the torus group \cite{bowditch1998plms} as well
as to the case of hyperbolic once punctured torus bundles
\cite{bowditch1997t}. There are also some other generalizations
along these directions, by Makoto Sakuma and his co-workers, see
\cite{akiyoshi-miyachi-sakuma2002preprint}, \cite{sakuma1999sk}.

\vskip 10pt

In this paper we further generalize McShane's identity as in
Theorem \ref{thm:mcshane general} to the cases of hyperbolic
cone-surfaces possibly with cusps and/or geodesic boundary. (See
for example \cite{cooper-hodgson-kerckhoff2000book} for basic
facts on cone-manifolds.) We assume that all cone points have cone
angle $\le \pi$ (except for the one-cone torus where we allow the
cone angle up to $ 2 \pi$). The ideas are related in spirit to
those in \cite{basmajian1993ajm} while the method of proof follows
closely that of McShane's in \cite{mcshane1998im}. The key points
are that the assumption that all cone angles are $\le \pi$ implies
that all non-peripheral simple closed curves are essentially
realizable as simple geodesics in their free (relative) homotopy
classes; and that the Birman-Series result
\cite{birman-series1985t} on the sparsity of simple geodesics
carries over to this case, in particular to simple geodesic rays
emanating (normally) from a fixed boundary component.  It should
be noted that our result shows that the assumption of discreteness
of the holonomy group is unnecessary, and that it gives identities
for all hyperbolic orbifold surfaces. We also show how the result
can be formulated in terms of complex lengths (Theorem
\ref{thm:complexified}) even though the situation we consider here
is real. This is particularly useful, and is explored further in
\cite{tan-wong-zhang2004preprint}, where we show how this approach
allows us to generalize McShane's identity to Schottky groups, and
how the Markoff triples and analytic continuation methods adopted
by Bowditch in \cite{bowditch1996blms} can be generalized as well.
(See also \cite{goldmanGT2003} for related work on generalized
Markoff triples.) This should also lead to generalizations of
Bowditch's interpretation \cite{bowditch1997t} of McShane's
identity for complete hyperbolic 3-manifolds which are once
punctured torus bundles over the circle to identities for the
hyperbolic 3-manifolds obtained by hyperbolic Dehn surgery on such
manifolds. This will be explored in future work, and should tie up
nicely with the work of Sakuma in \cite{sakuma1999sk}, and
Akiyoshi-Miyachi-Sakuma in
\cite{akiyoshi-miyachi-sakuma2002preprint} and
\cite{akiyoshi-miyachi-sakuma2004preprint}.

\vskip 10pt

To state the most general form of our generalized McShane's
identities, we need to introduce some new terminology. However, to
let the reader get the flavor of the generalized identities, we
first state the corresponding generalizations of Theorems
\ref{thm:mcshane torus} and \ref{thm:mcshane general}.

\vskip 10pt

\begin{thm}\label{thm:mcshane cone hole torus}
Let $T$ be either a hyperbolic one-cone torus where the single
cone point has cone angle $\theta \in (0, 2\pi)$ or a hyperbolic
one-hole torus where the single boundary geodesic has length
$l>0$. Then we have respectively
\begin{eqnarray}\label{eqn:mcshane cone torus}
\sum_{\gamma} 2 \tan^{-1} \left ( \frac {\sin \frac{\theta}{2}}
{\cos \frac{\theta}{2}+ \exp |\gamma|} \right )= \frac{\theta}{2},
\end{eqnarray}
\begin{eqnarray}
\sum_{\gamma} 2 \tanh^{-1}\left ( \frac {\sinh \frac{l}{2}} {\cosh
\frac{l}{2}+ \exp |\gamma|} \right )= \frac{l}{2},
\end{eqnarray}

\noindent where the sum in either case extends over all simple
closed geodesics on $T$.
\end{thm}

\vskip 10pt

\begin{thm}\label{thm:mcshane cone hole surface}
Let $M$ be a compact hyperbolic cone-surface with a single cone
point of cone angle $\theta \in (0, \pi]$ and without boundary or
let $M$ be a compact hyperbolic surface with a single boundary
geodesic having length $l>0$. Then we have respectively
\begin{eqnarray}
\sum 2 \tan^{-1} \left ( \frac {\sin
\frac{\theta}{2}}{\cos\frac{\theta}{2}+\exp
\frac{|\alpha|+|\beta|}{2}} \right )= \frac{\theta}{2},
\end{eqnarray}
\begin{eqnarray}
\sum 2 \tanh^{-1} \left ( \frac {\sinh \frac{l}{2}}{\cosh
\frac{l}{2}+ \exp \frac{|\alpha|+|\beta|}{2}} \right )=
\frac{l}{2},
\end{eqnarray}

\noindent where the sum in either case extends over all unordered
pairs of simple closed geodesics on $M$ which bound with the cone
point (respectively, the boundary geodesic) an embedded pair of
pants.
\end{thm}

\vskip 10pt

For the purposes of this paper we make the following definition.

\vskip 10pt

\begin{defn}
A {\bf compact hyperbolic cone-surface} $M$ is a compact
(topological) surface $M$ with hyperbolic cone structure where
each boundary component is a smooth simple closed geodesic and
where there are a finite number of interior points which form all
the cone points and cusps. Its {\bf geometric boundary}, denoted
$\Delta M$, is the union of all cusps, cone points and geodesic
boundary components. (Note that $\Delta M$ is different from the
usual topological boundary $\partial M$ when there are cusps or
cone points.) Thus a {\bf geometric boundary component} is either
a cusp, a cone point, or a boundary geodesic. The {\bf geometric
interior} of $M$ is $M - \Delta M$.
\end{defn}

\vskip 10pt

In this paper we consider a compact hyperbolic cone-surface $M =
M(\Delta_{0}; k, \Theta, L)$ with $k$ cusps $C_{1}, C_{2}, \cdots,
C_{k}$, with $m$ cone points $P_{1}, P_{2}, \cdots, P_{m}$, where
the cone angle of $P_{i}$ is $\theta_{i} \in (0, \pi]$, $i=1, 2,
\cdots, m$, and with $n$ geodesic boundary components $B_{1},
B_{2}, \cdots, B_{n}$, where the length of $B_{i}$ is $l_{i}>0$,
$i=1, 2, \cdots, n$, together with an extra {\it distinguished}
geometric boundary component $\Delta_{0}$. Thus $\Delta_{0}$ is
either a cusp $C_{0}$ or a cone point $P_{0}$ of cone angle
$\theta_{0} \in (0, \pi]$ or a geodesic boundary component $B_{0}$
of length $l_{0}>0$. Note that in the above notation $\Theta =
(\theta_{1}, \theta_{2}, \cdots, \theta_{m})$ and $L = (l_{1},
l_{2}, \cdots, l_{n})$. We exclude the case where $M$ is a
geometric pair of pants for we have only trivial identities in
that case.

\vskip 10pt

We allow that some (even all) of the cone angles $\theta_{i}$ are
equal to $\pi$, $i=0, 1, \cdots, m$; these are often cases of
particular interest. However, for clarity of exposition, quite
often in proofs/statements of lemmas/theorems we shall first
consider the case where all the cone angles are less than $\pi$
and then point out the addenda that should be made when there are
angle $\pi$ cone points. The advantage of this assumption of
strict inequality is that every non-trivial, non-peripheral simple
closed curve on such $M$ can be realized as a (smooth) simple
closed geodesic in its free homotopy class in the geometric
interior of $M$ under the given hyperbolic cone-structure (see \S
\ref{s:realizibility} for the proof of this statement).

\vskip 10pt

We call a simple closed curve on $M$ {\it peripheral} if it is
freely homotopic on $M$ to a geometric boundary component of $M$.

\vskip 10pt

\begin{defn}\label{defn:gscg}
By a {\bf generalized simple closed geodesic} on $M$ we mean
either \begin{itemize} \item[(i)] a simple closed geodesic in the
geometric interior of $M$; or \item[(ii)] a degenerate simple
closed geodesic which is the double of a simple geodesic arc in
the geometric interior of $M$ connecting two angle $\pi$ cone
points; or \item[(iii)] a geometric boundary component, that is, a
cusp or a cone point or a boundary geodesic. \end{itemize}
\end{defn}
\noindent In particular, generalized simple closed geodesics of
the first two kinds are called {\bf interior generalized simple
closed geodesics}.

\vskip 10pt

For each pair of generalized simple closed geodesics $\alpha,
\beta$ which bound with $\Delta_{0}$ an embedded geometric pair of
pants we shall define in \S \ref{s:dGf} a {\bf gap function} ${\rm
Gap}(\Delta_{0};\alpha, \beta)$ when $\Delta_{0}$ is a cone point
or a boundary geodesic as well as a {\bf normalized gap function}
${\rm Gap}^{\prime}(\Delta_{0};\alpha, \beta)$ when $\Delta_{0}$
is a cusp.

\vskip 10pt

Now we are in a position to state the most general (real) form of
our generalization of McShane's identity.

\vskip 10pt

\begin{thm}\label{thm:mcshane most general}
Let $M$ be a compact hyperbolic cone-surface with all cone angles
in $(0, \pi]$.  Then one has either
\begin{eqnarray}\label{eqn:001}
\sum {\rm Gap}(\Delta_{0};\alpha, \beta) = \frac{\theta_{0}}{2},
\end{eqnarray}
when $\Delta_{0}$ is a cone point of cone angle $\theta_{0}$;
\quad or
\begin{eqnarray}\label{eqn:002}
\sum {\rm Gap}(\Delta_{0};\alpha, \beta) = \frac{l_{0}}{2},
\end{eqnarray}
when $\Delta_{0}$ is a boundary geodesic of length $l_{0}$; \quad
or
\begin{eqnarray}\label{eqn:00'}
\sum {\rm Gap}^{\prime}(\Delta_{0};\alpha, \beta) = \frac{1}{2},
\end{eqnarray}
when $\Delta_{0}$ is a cusp; where in each case the sum is over
all pairs of generalized simple closed geodesics $\alpha, \beta$
on $M$ which bound with $\Delta_{0}$ an embedded pair of pants.
\end{thm}

\vskip 10pt

\vskip 10pt

\begin{rmk}~\\
\begin{itemize}
\item[(i)] In the case of the hyperbolic
one-cone torus, the theorem holds for $\theta_0 \in (0, 2\pi)$.
\item [(ii)] In the special cases where the geometric boundary
$\Delta M$ is a single cone point or a single boundary geodesic
Theorem \ref{thm:mcshane most general} gives all the previously
stated generalized identities in Theorems \ref{thm:mcshane cone
hole torus} and \ref{thm:mcshane cone hole surface}.
\item[(iii)] The cusp case (that is, $\Delta_{0}$ is a cusp) is the
limit case of the other cases as the cone angle $\theta_{0}$ or
the boundary geodesic length $l_{0}$ approaches $0$, and the
identity in the cusp case can indeed be derived from the first
order infinitesimal of the identities of the other cases.
\end{itemize}
\end{rmk}

\vskip 10pt

It is also interesting to note that McShane's Weierstrass
identities can be deduced as special cases of our general Theorem
\ref{thm:mcshane most general} by applying the theorem to the
quotient of the once punctured torus by its elliptic involution
and then lifting back to the torus. Thus we have the following
generalized Weierstrass identities:

\vskip 10pt

\begin{cor}\label{cor:mcshane conical holed weierstrass}
Let $T$ be either a hyperbolic one-cone torus where the single
cone point has cone angle $\theta \in (0, 2\pi)$ or a hyperbolic
one-hole torus where the single boundary geodesic has length
$l>0$. Then we have respectively
\begin{eqnarray}
\sum_{\gamma \in \mathcal A} \tan^{-1} \left ( \frac{\cos
\frac{\theta}{4}}{\sinh \frac{|\gamma|}{2}} \right )=
\frac{\pi}{2},
\end{eqnarray}
\begin{eqnarray}
\sum_{\gamma \in \mathcal A} \tan^{-1} \left ( \frac{\cosh
\frac{l}{4}}{\sinh \frac{|\gamma|}{2}} \right )= \frac{\pi}{2},
\end{eqnarray}
where the sum in either case is over all the simple closed
geodesics $\gamma$ in a Weierstrass class $\mathcal A$.
\end{cor}

\vskip 10pt

\noindent McShane's original Weierstrass identity
(\ref{eqn:mcshane weierstrass}) then corresponds to the case
$\theta = 0$ or $l = 0$ in the above two identities, noticing that
\begin{eqnarray*}
\tan^{-1} \left ( \frac{1}{\sinh \frac{|\gamma|}{2}}\right )=
\sin^{-1} \left ( \frac{1}{\cosh \frac{|\gamma|}{2}}\right ).
\end{eqnarray*}

\vskip 10pt

As further corollaries, there are the following weaker but neater
identities, each of which is obtained by summing the three
McShane's Weierstrass identities in the corresponding case.

\begin{cor}\label{cor:mcshane combined weierstrass}
Let $T$ be a hyperbolic torus whose geometric boundary is either a
single cusp, a single cone point of cone angle $\theta \in (0,
2\pi)$, or a single boundary geodesic of length $l>0$. Then we
have respectively
\begin{eqnarray}
\sum_{\gamma} \tan^{-1} \left ( \frac{1}{\sinh \frac{|\gamma|}{2}}
\right )= \frac{3\pi}{2},
\end{eqnarray}
\begin{eqnarray}
\sum_{\gamma} \tan^{-1} \left ( \frac{\cos \frac{\theta}{4}}{\sinh
\frac{|\gamma|}{2}}\right )= \frac{3\pi}{2},
\end{eqnarray}
\begin{eqnarray}\label{eqn:combinedWeierstrassonehole}
\sum_{\gamma} \tan^{-1} \left ( \frac{\cosh \frac{l}{4}}{\sinh
\frac{|\gamma|}{2}}\right )= \frac{3\pi}{2},
\end{eqnarray}
where the sum in each case is over all the simple closed geodesics
$\gamma$ on $T$.
\end{cor}

\vskip 10pt

\begin{rmk}
The identity (\ref{eqn:combinedWeierstrassonehole}) was also
obtained by McShane \cite{mcshane2004preprint} using Wolpert's
variation of length method. It seems likely his method can be
extended to prove some of the other identities as well.
\end{rmk}

Similarly, for a genus two closed hyperbolic surface $M$, one can
consider the (six) identities on the quotient  surface $M/\eta$
where $\eta$ is the unique hyper-elliptic involution on $M$ (note
that $M/\eta$ is a closed hyperbolic orbifold of genus $0$ with
six cone angle $\pi$ points, and we may choose any one of these
cone points to be the distinguished geometric boundary component)
and re-interpret them as Weierstrass identities on the original
surface $M$ (see also McShane \cite{mcshane1998preprint} where the
Weierstrass identities were obtained directly). Combining all the
six Weierstrass identities for $M$, we then have the following
very neat identity.

\vskip 10pt

\begin{thm}\label{thm:mcshane genus two global}
Let $M$ be a genus two closed hyperbolic surface. Then
\begin{eqnarray}
\sum \tan^{-1} \exp \left( -\frac{|\alpha|}{4}- \frac{|\beta|}{2}
\right ) = \frac{3\pi}{2},
\end{eqnarray}
where the sum is over all ordered pairs $(\alpha, \beta)$ of
disjoint simple closed geodesics on $M$ such that $\alpha$ is
separating and $\beta$ is non-separating.
\end{thm}

\vskip 10pt

\begin{rmk}
This is the only case that we know of where McShane's identity
extends in a nice way to a closed surface.
\end{rmk}

\vskip 10pt

We observe that the above identity for closed genus two surface
$M$ also extends to quasi-Fuchsian representations of
${\pi}_{1}(M)$. More precisely, let $\rho: {\pi}_{1}(M)
\rightarrow {\rm SL}(2, \mathbf C)$ be a quasi-Fuchsian
representation, that is, $\pi \circ \rho: {\pi}_{1}(M) \rightarrow
{\rm PSL}(2, \mathbf C)$ is a quasi-Fuchsian representation where
$\pi: {\rm SL}(2, \mathbf C) \rightarrow {\rm PSL}(2, \mathbf C)$
is the projection map. For each essential simple closed curve
$\gamma$, let $l_{\rho}(\gamma)/2 \in \mathbf C$ with positive
real part and with imaginary part $\in (-\pi, \pi]$ be defined by
$${\rm tr} \rho ([\gamma]) = 2 \cosh (l_{\rho}(\gamma)/2),$$ where
$[\gamma] \in {\pi}_{1}(M)$ is the homotopy class of $\gamma$.
Note that $l_{\rho}(\gamma)$ is also called the complex length of
$\rho([\gamma])$, see for example \cite{fenchel1989book}.

\vskip 10pt

\begin{add}\label{add:genus two}
For a quasi-Fuchsian representation $\rho: {\pi}_{1}(M)
\rightarrow {\rm SL}(2, \mathbf C)$ for the  closed genus two
surface $M$, we have
\begin{eqnarray}\label{eqn:genus two QF}
\sum \tan^{-1} \exp \left( -\frac{l_{\rho}(\alpha)}{4}-
\frac{l_{\rho}(\beta)}{2} \right ) = \frac{3\pi}{2},
\end{eqnarray}
where the sum is over all the ordered pairs $[\alpha],[\beta]$ of
homotopy classes of disjoint unoriented essential simple closed
curves $\alpha, \beta$ on $M$ such that $\alpha$ is non-separating
and $\beta$ is separating.
\end{add}

\vskip 15pt

In the statement of Theorem \ref{thm:mcshane most general} we did
not write down the explicit expression for the gap functions due
to their ``case by case'' nature as can be seen in \S \ref{s:dGf}.
The cone points and boundary geodesics as geometric boundary
components seem to have different roles in the series in the
generalized identities, hence making the identities not in a
unified form. This difference can, however, be removed by
assigning purely imaginary length to a cone point as a geometric
boundary component. More precisely, for each generalized simple
closed geodesic $\delta$, we define its {\bf complex length}
$|\delta|$ as: $|\delta|=0$ if $\delta$ is a cusp;
$|\delta|=\theta i$ if $\delta$ is a cone point of angle $\theta
\in (0, \pi]$; and $|\delta|=l$ if $\delta$ is a boundary geodesic
or an interior generalized simple closed geodesic of length $l>0$.
Then we can reformulate the generalized McShane's identities in
Theorem \ref{thm:mcshane most general} as follows.

\vskip 10pt

\begin{thm}\label{thm:complexified}
Let $M$ be a compact hyperbolic cone-surface with all cone angles
in $(0, \pi]$, and let all its geometric boundary components be
$\Delta_0, \Delta_1, \cdots, \Delta_N$ with complex lengths $L_0,
L_1, \cdots, L_N$ respectively. Then
\begin{eqnarray}\label{eqn:reform of cp and gb cases}
\lefteqn{\sum_{\alpha, \beta} 2 \tanh^{-1} \left ( \frac {\sinh
\frac{L_0}{2}}{\cosh \frac{L_0}{2}+
\exp\frac{|\alpha|+|\beta|}{2}} \right )}  \\ & & +
\sum_{j=1}^{N}\sum_{\beta} \tanh^{-1} \left (
\frac{\sinh\frac{L_0}{2}\sinh\frac{L_j}{2}
}{\cosh\frac{|\beta|}{2}+\cosh\frac{L_0}{2}\cosh\frac{L_j}{2} }
\right ) = \frac{L_0}{2}\nonumber,
\end{eqnarray}
if $\Delta_0$ is a cone point or a boundary geodesic; and
\begin{eqnarray}\label{eqn:reform of cusp cases}
\sum_{\alpha,\beta}\frac{1}{1+\exp\frac{|\alpha|+|\beta|}{2}}+
\sum_{j=1}^{N}\sum_{\beta}\frac{1}{2}\frac{\sinh\frac{L_j}{2}}
{\cosh\frac{|\beta|}{2}+\cosh\frac{L_j}{2}}=\frac{1}{2},
\end{eqnarray}
if $\Delta_0$ is a cusp; where in either case the first sum is
over all (unordered) pairs of generalized simple closed geodesics
$\alpha, \beta$ on $M$ which bound with $\Delta_0$ an embedded
pair of pants on $M$ (note that one of $\alpha, \beta$ might be a
geometric boundary component) and the sub-sum in the second sum is
over all interior simple closed geodesics $\beta$ which bounds
with $\Delta_j$ and $\Delta_0$ an embedded pair of pants on $M$.

Furthermore, each series in (\ref{eqn:reform of cp and gb cases})
and (\ref{eqn:reform of cusp cases}) converges absolutely.
\end{thm}

\bigskip

\noindent {\it Additional Remark.} We were informed while writing
this paper by Makoto Sakuma and Caroline Series of the recent
striking results of Maryam Mirzakhani
\cite{mirzakhani2004preprint} where she had generalized McShane's
identities to hyperbolic surfaces with boundary and used it to
calculate the Weil-Petersson volumes of the corresponding moduli
spaces. There is obviously an overlap of her results with ours, in
particular, the identities she obtains are equivalent to ours in
the case of hyperbolic surfaces with boundary (see \S
\ref{s:reformulation} for further explanations). In fact, her
expressions in terms of the $\log$ function seems particular well
suited to her purpose of calculating the Weil-Petersson volumes.
It also seems (as already observed by her in
\cite{mirzakhani2004preprint}) that her methods should extend
fairly easily to cover the case of volumes of the moduli spaces of
compact hyperbolic cone-surfaces with all cone angles bounded
above by $\pi$, as defined and used in our context, and that the
formulas she exhibited for the volumes should hold in this case as
well, using the convention that a cone point of angle $\theta$
corresponds to a geometric boundary component with purely
imaginary length $\theta i$ .

\bigskip
\noindent {\it Acknowledgements.} We would like to thank Caroline
Series for helpful conversations; Makoto Sakuma for his
encouragement to write up our results on the cone-manifold case
(during conversations with the first named author at the Isaac
Newton Institute in Aug, 2003) and also for bringing to our
attention the recent works of McShane \cite{mcshane2004preprint}
and Mirzakhani \cite{mirzakhani2004preprint}; and Greg McShane for
helpful e-mail correspondence and also for bringing our attention
to \cite{mcshane1998preprint}.

\vskip 30pt
\section{{\bf The organization of the rest of this paper}}
\vskip 30pt

The rest of this paper is organized as follows. In \S \ref{s:dGf}
we define the gap functions used in Theorem \ref{thm:mcshane most
general} for the various cases. In \S \ref{s:realizibility} we
deal with the problem of realization of simple closed curves by
geodesics, and show that the assumption that all cone angles are
less than or equal to $\pi$ is essential. In \S \ref{s:gaps} we
analyze the so-called $\Delta_{0}$-geodesics, that is, the
geodesics starting/emanating orthogonally from $\Delta_{0}$, and
determine all the gaps between all simple-normal
$\Delta_{0}$-geodesics. In \S \ref{s:calculation} we calculate the
gap function which is the width of a combined gap measured
suitably. In \S \ref{s:gBS} we generalize the Birman-Series
theorem (which states that the point set of all complete geodesics
with bounded self intersection numbers on a compact hyperbolic
surface has Hausdorff dimension 1) to the case of compact
hyperbolic cone-surfaces with all cone angles less than or equal
to $\pi$. We prove the theorems in this paper in \S \ref{s:proof},
except for Theorem \ref{thm:complexified}, which is deferred to
the last section. Finally in \S \ref{s:reformulation} we restate
the complexified generalized McShane's identity (\ref{eqn:reform
of cp and gb cases}) (Theorem \ref{thm:complexified}) using two
functions of complex variables and hence unify the somewhat
unattractive ``case-by-case'' definition of the gap functions. We
interpret the geometric meanings of the complexified summands in
the complexified generalized McShane's identity and prove the
absolute convergence of the complexified series in it by a simple
use of the Birman--Series arguments in \cite{birman-series1985t}.

\vskip 20pt

\vskip 30pt
\section{{\bf Defining the Gap functions}}\label{s:dGf}
\vskip 30pt

In this section, for a compact hyperbolic cone-surface $M =
M(\Delta_{0}; k, \Theta, L)$ with all cone angles $\le \pi$ we
define the gap function ${\rm Gap}(\Delta_{0};\alpha, \beta)$
(when $\Delta_{0}$ is a cone point or a boundary geodesic) and the
normalized gap function ${\rm Gap}^{\prime}(\Delta_{0};\alpha,
\beta)$ (when $\Delta_{0}$ is a cusp) where $\alpha, \beta$ are
generalized simple closed geodesics on $M$ which bound with
$\Delta_{0}$ a geometric pair of pants.

\vskip 10pt

Throughout this paper we use $|\alpha|$ to denote the length of
$\alpha$ when $\alpha$ is an interior generalized simple closed
geodesic or a boundary geodesic. In particular, when $\alpha$ is a
degenerate simple closed geodesic (that is, the double cover of a
simple geodesic arc which connects two angle $\pi$ cone points),
its length $|\alpha|$ is defined as twice the length of the simple
geodesic that it covers.

\vskip 10pt

Recall that an interior generalized simple closed geodesic is
either a simple closed geodesic in the geometric interior of $M$
or a degenerate simple closed geodesic on $M$ which is the double
cover of a simple geodesic arc which connects two angle $\pi$ cone
points.

\vskip 20pt

{\it Case} 0. $\Delta_{0}$ is a cusp.

\vskip 10pt

{\it Subcase} 0.1.  Both $\alpha$ and $\beta$ are interior
generalized simple closed geodesics.

\vskip 10pt

In this case
\begin{eqnarray}
{\rm Gap}^{\prime}(\Delta_{0};\alpha, \beta)=
\frac{1}{1+\exp\frac{1}{2}(|\alpha|+|\beta|)}.
\end{eqnarray}

\vskip 10pt

{\it Subcase} 0.2. One of $\alpha, \beta$, say $\alpha$, is a
boundary geodesic and the other, $\beta$, is an interior
generalized simple closed geodesic.

\vskip 10pt

In this case
\begin{eqnarray}
{\rm Gap}^{\prime}(\Delta_{0};\alpha, \beta)=
\frac{1}{2}-\frac{1}{2}
\frac{\sinh\frac{|\beta|}{2}}{\cosh\frac{|\alpha|}{2}+\cosh\frac{|\beta|}{2}}.
\end{eqnarray}

\vskip 10pt

{\it Subcase} 0.3. One of $\alpha, \beta$, say $\alpha$, is a cone
point of cone angle $\varphi \in (0,\pi]$ and the other, $\beta$,
is an interior generalized simple closed geodesic.

\vskip 10pt

In this case
\begin{eqnarray}
{\rm Gap}^{\prime}(\Delta_{0};\alpha, \beta)=
\frac{1}{2}-\frac{1}{2}
\frac{\sinh\frac{|\beta|}{2}}{\cos\frac{\varphi}{2}+\cosh\frac{|\beta|}{2}}.
\end{eqnarray}

\vskip 10pt

{\it Subcase} 0.4. One of $\alpha, \beta$, say $\alpha$, is also a
cusp and the other, $\beta$, is an interior generalized simple
closed geodesic.

\vskip 10pt

In this case
\begin{eqnarray}
{\rm Gap}^{\prime}(\Delta_{0};\alpha, \beta) =
\frac{1}{2}-\frac{1}{2}
\frac{\sinh\frac{|\beta|}{2}}{1+\cosh\frac{|\beta|}{2}} =
\frac{1}{1+\exp \frac{1}{2}|\beta|},
\end{eqnarray}
which is the common value of ${\rm Gap}(\Delta_{0};\alpha, \beta)$
in Subcases 0.1 through 0.3 when $|\alpha|=0$.

\vskip 20pt

{\it Case} 1. $\Delta_{0}$ is a cone point of cone angle $\theta
\in (0,\pi]$.

\vskip 10pt

{\it Subcase} 1.1. Both $\alpha$ and $\beta$ are interior
generalized simple closed geodesics.

\vskip 10pt

In this case
\begin{eqnarray}\label{eqn:dGf 1.1}
{\rm Gap}(\Delta_{0};\alpha, \beta) = 2 \tan^{-1} \left
(\frac{\sin\frac{\theta}{2}}{\cos\frac{\theta}{2}+
\exp\frac{|\alpha|+|\beta|}{2}} \right ).
\end{eqnarray}

\vskip 10pt

{\it Subcase} 1.2. One of $\alpha, \beta$, say $\alpha$, is a
boundary geodesic and the other, $\beta$, is an interior
generalized simple closed geodesic.

\vskip 10pt

In this case
\begin{eqnarray}
{\rm Gap}(\Delta_{0};\alpha, \beta)= \frac{\theta}{2} - \tan^{-1}
\left ( \frac{\sin\frac{\theta}{2}\sinh\frac{|\beta|}{2}
}{\cosh\frac{|\alpha|}{2}+\cos\frac{\theta}{2}\cosh\frac{|\beta|}{2}
} \right ).
\end{eqnarray}

\vskip 10pt

{\it Subcase} 1.3. One of $\alpha, \beta$, say $\alpha$, is a cone
point of cone angle $\varphi \in (0,\pi]$ and the other, $\beta$,
is an interior generalized simple closed geodesic.

\vskip 10pt

In this case
\begin{eqnarray}\label{subcase 1.3}
{\rm Gap}(\Delta_{0};\alpha, \beta)= \frac{\theta}{2} - \tan^{-1}
\left ( \frac{\sin\frac{\theta}{2}\sinh\frac{|\beta|}{2}
}{\cos\frac{\varphi}{2}+\cos\frac{\theta}{2}\cosh\frac{|\beta|}{2}
} \right ).
\end{eqnarray}

\vskip 10pt

Note that there is no gap when $\theta=\varphi=\pi$.

\vskip 10pt

{\it Subcase} 1.4. One of $\alpha, \beta$, say $\alpha$, is a cusp
and the other, $\beta$, is an interior generalized simple closed
geodesic.

\vskip 10pt

In this case
\begin{eqnarray}
{\rm Gap}(\Delta_{0};\alpha, \beta) &=& 2 \tan^{-1} \left ( \frac
{\sin \frac{\theta}{2}}{\cos\frac{\theta}{2}+
\exp\frac{|\beta|}{2}}
\right )\\
&=& \frac{\theta}{2} - \tan^{-1} \left (
\frac{\sin\frac{\theta}{2}\sinh\frac{|\beta|}{2}
}{1+\cos\frac{\theta}{2}\cosh\frac{|\beta|}{2} } \right ),
\end{eqnarray}
which is the common value of ${\rm Gap}(\Delta_{0};\alpha, \beta)$
in Subcases 1.1 through 1.3 when $|\alpha|=0$.

\vskip 20pt

{\it Case} 2. $\Delta_{0}$ is a boundary geodesic of length $l>0$.

\vskip 10pt

{\it Subcase} 2.1. Both $\alpha$ and $\beta$ are interior
generalized simple closed geodesics.

\vskip 10pt

In this case
\begin{eqnarray}
{\rm Gap}(\Delta_{0};\alpha, \beta)= 2 \tanh^{-1} \left ( \frac
{\sinh \frac{l}{2}}{\cosh \frac{l}{2}+
\exp\frac{|\alpha|+|\beta|}{2}} \right ).
\end{eqnarray}

\vskip 10pt

{\it Subcase} 2.2. One of $\alpha, \beta$, say $\alpha$, is a
boundary geodesic and the other, $\beta$, is an interior
generalized simple closed geodesic.

\vskip 10pt

In this case
\begin{eqnarray}
{\rm Gap}(\Delta_{0};\alpha, \beta)= \frac{l}{2} - \tanh^{-1}\left
( \frac{\sinh\frac{l}{2}\sinh\frac{|\beta|}{2}
}{\cosh\frac{|\alpha|}{2}+\cosh\frac{l}{2}\cosh\frac{|\beta|}{2}
}\right ).
\end{eqnarray}

\vskip 10pt

{\it Subcase} 2.3. One of $\alpha, \beta$, say $\alpha$, is a cone
point of cone angle $\varphi \in (0,\pi]$ and the other, $\beta$,
is an interior generalized simple closed geodesic.

\vskip 10pt

In this case
\begin{eqnarray}
{\rm Gap}(\Delta_{0};\alpha, \beta)= \frac{l}{2} - \tanh^{-1}
\left ( \frac{\sinh\frac{l}{2}\sinh\frac{|\beta|}{2}
}{\cos\frac{\varphi}{2}+\cosh\frac{l}{2}\cosh\frac{|\beta|}{2}
}\right ).
\end{eqnarray}

\vskip 10pt

{\it Subcase} 2.4. One of $\alpha, \beta$, say $\alpha$, is a cusp
and the other, $\beta$, is an interior generalized simple closed
geodesic.

\vskip 10pt

In this case
\begin{eqnarray}
{\rm Gap}(\Delta_{0};\alpha, \beta) &=& 2 \tanh^{-1} \left ( \frac
{\sinh \frac{l}{2}}{\cosh \frac{l}{2}+ \exp \frac{|\beta|}{2}}
\right )\\
&=& \frac{l}{2} - \tanh^{-1}\left (
\frac{\sinh\frac{l}{2}\sinh\frac{|\beta|}{2}
}{1+\cosh\frac{l}{2}\cosh\frac{|\beta|}{2} }\right ),
\end{eqnarray}
which is the common value of ${\rm Gap}(\Delta_{0};\alpha, \beta)$
in Subcases 2.1 through 2.3 when $|\alpha|=0$.

\vskip 30pt

\vskip 30pt
\section{{\bf Realizing simple curves by geodesics on hyperbolic
cone-surfaces}}\label{s:realizibility}
\vskip 30pt

In this section we consider the problem of realizing essential
simple curves in their free (relative) homotopy classes by
geodesics on a compact hyperbolic cone-surface $M$ with all cone
angles smaller than $\pi$. We show that each essential simple
closed curve in the geometric interior of $M$ can be realized
uniquely in its free homotopy class (where the homotopy takes
place in the geometric interior of $M$) as either a geometric
boundary component or a simple closed geodesic in the geometric
interior of $M$. We also show that each essential simple arc which
connects geometric boundary components of $M$ can be realized
uniquely in its free relative homotopy class (where the homotopy
takes place in the geometric interior of $M$ and the endpoints
slide on the same geometric boundary components) as a simple
geodesic arc which is normal to the geometric boundary components
involved. We also make addenda for the cases when there are angle
$\pi$ cone points.

\vskip 10pt

\begin{thm}\label{thm:realization}
Let $M$ be a compact hyperbolic cone-surface with all cone angles
less than $\pi$.

(i) If $c$ is an essential non-peripheral simple closed curve in
the geometric interior of $M$, then there is a unique simple
closed geodesic in the free homotopy class of $c$ in the geometric
interior of $M$.

(ii) If $c$ is an essential simple arc which connects geometric
boundary components, then there is a unique simple normal geodesic
arc in the free relative homotopy class of $c$ in the geometric
interior of $M$ with endpoints varying on the respective geometric
boundary components.
\end{thm}

\begin{add}\label{add:realization}
If in addition $M$ has some cone angles equal to $\pi$, then
\begin{itemize} \item[(i)] in Theorem \ref{thm:realization}(i), if the
simple closed curve $c$ bounds with two angle $\pi$ cone points an
embedded pair of pants, then the geodesic realization for $c$ is
the double cover of the simple geodesic arc which connects these
two angle $\pi$ cone points and is homotopic (relative to
boundary) to a simple arc lying wholly in the pair of pants;
 \item[(ii)] in Theorem \ref{thm:realization}(ii), if the simple arc $c$
connects a geometric boundary component $\Delta$ to itself and
bounds together with $\Delta$ and an angle $\pi$ cone point $P$ an
embedded cylinder then the geodesic realization for $c$ is the
double cover of the normal simple geodesic arc which connects
$\Delta$ to $P$ and is homotopic (relative to boundary) to a
simple arc lying wholly in the cylinder.
\end{itemize}
\end{add}

The simple geodesic in Theorem \ref{thm:realization} and Addendum
\ref{add:realization} is called the geodesic realization of the
given simple curve in the respective homotopy class.

\vskip 10pt

The proof is a well-known use of the Arzela-Ascoli Theorem as used
in \cite{buser1992book} with slight modifications.

\vskip 10pt

\begin{pf}
(i) Suppose $c$ is an essential non-peripheral simple closed curve
in the geometric interior of $M$, parameterized on $[0, 1]$ with
constant speed. Let the length of $c$ be $|c|>0$. Then for each
cusp $C_{i}$, there is an embedded neighborhood $N(C_{i})$ of
$C_{i}$ on $M$, bounded by a horocycle, such that each
non-peripheral simple closed curve $c^{\prime}$ in the geometric
interior of $M$ with length $\le |c|$ cannot enter $N(C_{i})$; for
otherwise $c^{\prime}$ would be either peripheral or of infinite
length. Now let $M_{0}$ be $M$ with all the chosen horocycle
neighborhoods $N(C_{i})$ removed. Then $M_{0}$ is a compact metric
subspace of $M$ with the induced hyperbolic metric. Now choose a
sequence of simple closed curves $\{c_{k}\}_{1}^{\infty}$, where
each $c_{k}$ is parameterized on $[0, 1]$ with constant speed, in
the free homotopy class of $c$ (where the homotopy takes place in
the geometric interior of $M$) such that their lengths $\le |c|$
and are decreasing with limit the infimum of the lengths of the
simple closed curves in the free homotopy class of $c$. Then by
the Arzela-Ascoli Theorem (c.f. \cite{buser1992book} Theorem A.19,
page 429) there is a subsequence of $\{c_{k}\}_{1}^{\infty}$,
assumed to be $\{c_{k}\}_{1}^{\infty}$ itself, such that it
converges uniformly to a closed curve $\gamma$ in $M_{0}$. It is
clear that $\gamma$ is a geodesic since it is locally minimizing.
Note that $\gamma$ is away from cusps by the choice of
$\{c_{k}\}_{1}^{\infty}$. We claim that $\gamma$ cannot pass
through any cone point. For otherwise, suppose $\gamma$ passes
through a cone point $P$. Then for sufficiently large $k$, $c_{k}$
can be modified in the free homotopy class of $c$ to have length
smaller than $|\gamma|$ (since the cone point has cone angle
smaller than $\pi$), which is a contradiction. Thus $\gamma$ must
be a  closed geodesic in the geometric interior of $M$. The
uniqueness and simplicity of $\gamma$ can be proved by an easy
argument since there are no bi-gons in the hyperbolic plane.

\vskip 10pt

(ii) For an essential simple arc $c$ in the geometric interior of
$M$ which connects geometric boundary components, the proof of
case (i) applies without modifications when none of the involved
geometric boundary components is a cusp. Now suppose at least one
of the involved geometric boundary components is a cusp. For
definiteness let us assume that $c$ connects cusps $C_{1}$ to
$C_{2}$. Remove suitable horocycle neighborhoods $N(C_{1})$ and
$N(C_{2})$ respectively for $C_{1}$ and $C_{2}$ where the two
horocycles are $H_{1}$ and $H_{2}$ respectively. Choose a simple
arc $c_{0}$ in $M - N(C_{1}) \cup N(C_{2})$ which goes along $c$
and connects $H_{1}$ to $H_{2}$. Let the length of $c_{0}$ be
$|c_{0}|>0$. Now for all other cusps $C_{i}$, there is a horocycle
neighborhood $N(C_{i})$ of $C_{i}$ on $M$ such that each
non-peripheral simple closed curve $c^{\prime}$ in the geometric
interior of $M$ with length $\le |c_{0}|$ cannot enter $N(C_{i})$.
Again let $M_{0}$ be $M$ with all the chosen horocycle
neighborhoods $N(C_{i})$ removed. By the same argument as in (i)
we have a shortest simple geodesic realization $\gamma_{0}$ in the
free relative homotopy class of $c_{0}$ in $M_{0}$ and $c_{0}$
does not pass through any cone point. Hence $\gamma_{0}$ must be
perpendicular to both $H_{1}$ and $H_{2}$ at its endpoints. Thus
$\gamma_{0}$ can be extended to a  geodesic arc connecting $C_{1}$
to $C_{2}$. Again simplicity and uniqueness can be proved easily.
\end{pf}

\vskip 10pt

The addendum can be verified easily since the realizations as
degenerate simple geodesics in the respective cases are already
known.

\vskip 10pt

\begin{rmk}
We make a remark that the following fact, whose proof is easy and
hence omitted, is implicitly used through out this paper: On a
hyperbolic cone-surface for each cone point $P$ with angle less
than $\pi$ there is a cone region $N(P)$, bounded by a suitable
circle centered at $P$, such that if a geodesic $\gamma$ goes into
$N(P)$ then either $\gamma$ will go directly to the cone point $P$
(hence perpendicular to all the circles centered at $P$) or
$\gamma$ will develop a self-intersection in $N(P)$. The analogous
fact for a cusp is used in \cite{birman-series1985t},
\cite{haas1986actam} and \cite{mcshane1998im}.
\end{rmk}

\vskip 10pt

\vskip 30pt
\section{{\bf Gaps between simple-normal
$\Delta_{0}$-geodesics}}\label{s:gaps} \vskip 30pt

\begin{defn}\label{defn:deltageodesic}
A {\bf $\Delta_{0}$-geodesic} on $M$ is an {\it oriented} geodesic
ray which starts from $\Delta_{0}$ (and is perpendicular to it if
$\Delta_{0}$ is a boundary geodesic) and is fully developed, that
is, it develops forever until it terminates at a geometric
boundary component. We denote by ${\mathcal G}(\Delta_0)$ (or just
${\mathcal G}$) the set of $\Delta_0$-geodesics.
\end{defn}


A $\Delta_{0}$-geodesic is either non-simple or simple. It is
regarded as {\bf non-simple} if and only if it intersects itself
transversely at an interior point (a cone point is not treated as
an interior point) or at a point on a boundary geodesic. We shall
see later that somewhat surprisingly, in some sense, the set of
non-simple $\Delta_{0}$-geodesics is  easier to analyze than the
set of simple $\Delta_{0}$-geodesics.

\vskip 10pt

A simple $\Delta_{0}$-geodesic is either normal or not-normal in
the following sense:

A simple $\Delta_{0}$-geodesic is {\bf normal} if when fully
developed either it never intersects any boundary geodesic or it
intersects (hence terminates at) a boundary geodesic
perpendicularly. Note that a simple-normal $\Delta_{0}$-geodesic
may terminate at a cusp or a cone point. Thus a simple
$\Delta_{0}$-geodesic is {\bf not-normal} if and only if it
intersects a boundary geodesic (which might be $\Delta_{0}$
itself) obliquely.

\vskip 10pt

We shall analyze the structure of all non-simple and
simple-not-normal $\Delta_{0}$-geodesics and show that they form
gaps between simple-normal $\Delta_{0}$-geodesics. Furthermore,
the naturally measured widths of the suitably combined gaps are
given by the Gap functions defined before in \S \ref{s:dGf}.

\vskip 10pt

Note that McShane \cite{mcshane1998im} analyzes directly all
simple $\Delta_{0}$-geodesics (there are no simple-not-normal
$\Delta_{0}$-geodesics in his case since there are no geodesic
boundary componenets). Our analysis of the structure of
$\Delta_{0}$-geodesics is a bit different from and actually
simpler than that of McShane's. We shall analyze all non-simple
and simple-not-normal $\Delta_{0}$-geodesics and show that they
arise in the nice ways we expect.

\vskip 10pt

First we parameterize all the $\Delta_{0}$-geodesics and define
the widths for gaps between simple-normal $\Delta_{0}$-geodesics.

\vskip 8pt

If $\Delta_{0}$ is a cusp let $\mathcal H$ be a suitably chosen
small horocycle as in McShane \cite{mcshane1998im}, see also
\cite{haas1986actam}. If $\Delta_{0}$ is a cone point let
$\mathcal H$ be a suitably chosen small circle centered at
$\Delta_{0}$. Let $\mathcal H$ be $\Delta_{0}$ itself if
$\Delta_{0}$ is a boundary geodesic.

\vskip 8pt

Then each $\Delta_{0}$-geodesic has a unique first intersection
point with $\mathcal H$, which is the starting point when
$\Delta_{0}$ is a boundary geodesic. Note that the
$\Delta_{0}$-geodesics intersect $\mathcal H$ orthogonally at
their first intersection points. Thus ${\mathcal G}$ can be
naturally identified with ${\mathcal H}$, with the induced
topology and measure. Let ${\mathcal H}_{\texttt{ns}}$, ${\mathcal
H}_{\texttt{sn}}$, ${\mathcal H}_{\texttt{snn}}$ be the point sets
of the first intersections of $\mathcal H$ with respectively all
non-simple, all simple-normal, all simple-not-normal
$\Delta_{0}$-geodesics.

\begin{prop}
The set ${\mathcal H}_{\texttt{ns}}
\cup {\mathcal H}_{\texttt{snn}}$ is an open subset of $\mathcal
H$ and hence ${\mathcal H}_{\texttt{sn}}$ is a closed subset of
$\mathcal H$.
\end{prop}

\vskip 10pt \noindent {\it Proof}. It is easy to see that the
condition that either self-intersecting or ending obliquely at a
boundary component is an open condition. \square

\vskip 10pt

For the open subset ${\mathcal H}_{\texttt{ns}} \cup {\mathcal
H}_{\texttt{snn}}$ of $\mathcal H$, we determine its structure by
determining its maximal open intervals (which are the gaps we are
looking for). By a generalized Birman--Series Theorem (see \S
\ref{s:gBS}), the subset ${\mathcal H}_{\texttt{sn}}$ of $\mathcal
H$ has Hausdorff dimension $0$, and hence Lebesgue measure $0$.
Therefore the open subset ${\mathcal H}_{\texttt{ns}} \cup
{\mathcal H}_{\texttt{snn}}$ of $\mathcal H$ has full measure, and
our generalized McShane's identities
(\ref{eqn:001})-(\ref{eqn:00'}) follow immediately.

\vskip 10pt

\begin{defn}A $[\Delta_{0},\Delta_{0}]$-geodesic, $\gamma$, is an
(oriented) $\Delta_{0}$-geodesic which terminates at $\Delta_{0}$
perpendicularly. (With the orientation one can refer to its
starting point and ending point.) Hence the same geodesic with
reversed orientation (hence with the starting and ending points
interchanged) is also a $[\Delta_{0},\Delta_{0}]$-geodesic,
denoted by $-\gamma$.
\end{defn}

We say that a $[\Delta_{0},\Delta_{0}]$-geodesic $\gamma$ is a
{\bf degenerate simple} $[\Delta_{0},\Delta_{0}]$-{\bf geodesic}
if $\Delta_0$ is not a $\pi$ cone point, and $\gamma$ is the
double cover of a simple geodesic arc which connects $\Delta_{0}$
to an angle $\pi$ cone point, that is, $\gamma$ reaches the angle
$\pi$ cone point along the simple geodesic arc and goes back to
$\Delta_{0}$ along the same arc. Note that in this case
$\gamma=-\gamma$.

\vskip 10pt

We show that each non-degenerate simple
$[\Delta_{0},\Delta_{0}]$-geodesic $\gamma$ determines two maximal
open intervals of ${\mathcal H}_{\texttt{ns}} \cup {\mathcal
H}_{\texttt{snn}}$ as follows. (Their union is the {\it main gap},
defined later, determined by $\gamma$.)

\vskip 10pt

Consider the configuration $\gamma \cup \mathcal H$. Assume
$\gamma$ is non-degenerate and let $\mathcal H_1$ and $\mathcal
H_2$ be the two sub-arcs with endpoints inclusive that $\gamma$
divides $\mathcal H$ into. Note that $\gamma$ intersects $\mathcal
H$ twice (if ${\mathcal H}$ is taken to be a suitably small circle
about $\Delta_0$ when $\Delta_0$ is a cone point). Let $\gamma_0$
be the sub-arc of $\gamma$ between the two intersection points.
Thus we have two simple closed curves $\mathcal H_1 \cup \gamma_0$
and $\mathcal H_2 \cup \gamma_0$ on $M$. Their geodesic
realizations are disjoint generalized simple closed geodesics,
denoted $\alpha, \beta$ respectively (except when $M$ is a
hyperbolic torus with a single geometric boundary component, in
which case $\alpha=\beta$). Note that $\alpha, \beta$ bound with
$\Delta_{0}$ an embedded geometric pair of pants, denoted
$\mathcal P(\gamma)$, on $M$.

\vskip 10pt

Let $\delta_{\alpha}$ be the simple $\Delta_{0}$-geodesic arc in
$\mathcal P(\gamma)$ which terminates at $\alpha$ and is normal to
$\alpha$. Similarly, let $\delta_{\beta}$ be the simple
$\Delta_{0}$-geodesic arc in $\mathcal P(\gamma)$ which terminates
at $\beta$ and is normal to $\beta$. Let $[\alpha, \beta]$ be the
simple geodesic arc in $\mathcal P(\gamma)$ which connects
$\alpha$ and $\beta$ and is normal to them. See Figure
\ref{fig01}.

\vskip 5pt

\begin{figure}
\begin{center}
\mbox{\beginpicture \setcoordinatesystem units <0.25in,0.25in>
\setplotarea x from 1 to 9, y from 0 to 6 \ellipticalarc axes
ratio 2:1 360 degrees from 3.6 5 center at 5 5 \setquadratic
\setdashes<2pt> \plot 9.560660172  2.560660172  9.344866759
2.672628777  9.046371901 2.669812307  8.694394354  2.552486458
8.323388134  2.332135902 7.969669914  2.030330086  7.667864098
1.676611866  7.447513542 1.305605646  7.330187693  .9536280995
7.327371223  .6551332412 7.439339828  .439339828 / \setsolid \plot
9.560660172  2.560660172  9.672628777  2.344866759  9.669812307
2.046371901  9.552486458  1.694394354  9.332135902  1.323388134
9.030330086  .9696699142  8.676611866  .6678640980  8.305605646
.4475135417  7.953628099  .3301876929  7.655133241  .3273712228
7.439339828  .439339828 / \plot 2.560660172  .439339828
2.669812307  .9536280995  2.332135902 1.676611866  1.676611866
2.332135902  .9536280995  2.669812307 .439339828  2.560660172
.3301876929  2.046371901  .6678640980 1.323388134  1.323388134
.6678640980  2.046371901  .3301876929 2.560660172  .439339828 /
\plot 9.560660172  2.560660172 7.5 4 6.4 5 / \plot .439339828
2.560660172 2.5 4 3.6 5 / \plot 2.560660172  .439339828 5 1
7.439339828  .439339828 /
\plot 5.224157199 4.3 5.1 1.7 5 1  / \setdashes<2pt> \plot 5 1 4.8
2  4.775842801 5.7 / \setlinear\setsolid \plot 4.628894061 5.66
4.62 5.5 4.76 5.5 / \plot 5.39 4.3 5.39 4.15 5.21 4.14 /
\setquadratic \plot 5.75 4.4 5.9 1.7 6 0.9 / \setlinear \plot 5.9
4.45 5.9 4.3 5.8 4.27 / \setquadratic \setdashes<2pt> \plot 6 0.9
7 3 8 3.6 / \setsolid \plot 8 3.6 8.5 2.8 7.8 1.4 /
\plot 4.25 4.4 4.1 1.7 3.75 0.85 / \setdashes<2pt> \plot 3.75 0.85
2.4 2.1 2 3.6 / \setsolid \plot 2 3.6 3 2.6  3.5 0.78 /
\setdashes<2pt> \plot 3.5 0.78 3 1. 2.5 1.5 / \setlinear \setsolid
\plot 4.4 4.35 4.4 4.2 4.25 4.23 / \arrow <6pt> [.16,.6] from 5.19
3.7 to 5.188 3.6 \arrow <6pt> [.16,.6] from  4.23 3.7 to 4.22 3.6
\arrow <6pt> [.16,.6] from  5.79 3.7 to 5.8 3.6 \put {\mbox{\small
$\triangle_0$}} [cb] <0mm,2mm> at 5 5.5 \put {\mbox{\small
$\delta_\alpha$}} [rb] <-1mm,1mm> at 3 4.2 \put {\mbox{\small
$\delta_\beta$}} [lb] <1mm,1mm> at 7 4.2 \put {\mbox{\small
$\alpha$}} [rt] <-1mm,-1.6mm> at 1 1.2 \put {\mbox{\small
$\beta$}} [lt] <1mm,-1mm> at 9 1.2 \put {\mbox{\small $\gamma$}}
[lc] <1mm,0mm> at 5.1 3.2 \put {\mbox{\small
$\gamma_{{}_{\alpha}}$}} [rc] <-1mm,0mm> at 4.3 3.2 \put
{\mbox{\small $\gamma_{{}_{\beta}}$}} [lc] <1.5mm,0mm> at 5.7 3.2
\endpicture}
\hspace{0.5in} \mbox{\beginpicture \setcoordinatesystem units
<0.3in,0.3in> \setplotarea x from 3.5 to 9, y from 0.2 to 6 \plot
9.560660172 2.560660172 7 5 4 2  / \setquadratic \plot 9.560660172
2.560660172 9.639852575  2.377642961  9.607468267  2.108715941
9.466677248 1.780203564  9.231261125  1.424262911  8.924264069
1.075735931 8.575737089  .7687388748  8.219796436  .5333227521
7.891284059 .3925317336  7.622357039  .3601474247  7.439339828
.439339828 / \setdashes<2pt> \plot 9.560660172  2.560660172
9.377642961 2.639852575  9.108715941  2.607468267  8.780203564
2.466677248 8.424262911  2.231261125  8.075735931  1.924264069
7.768738875 1.575737089  7.533322752  1.219796435  7.392531733
.8912840588 7.360147425  .6223570393  7.439339828  .439339828 /
\put {\mbox{\LARGE $\cdot$}} [cc] <0mm,0mm> at 7 5 \put
{\mbox{\scriptsize $\triangle_{0}$}} [cb] <0mm,1mm> at 7 5 \put
{\mbox{\LARGE $\cdot$}} [cc] <0mm,0mm> at 4 2 \put
{\mbox{\scriptsize $\alpha$}} [rt] <0mm,0mm> at 4 2 \setsolid
\plot 4 2 6.3 1.2 7.439339828  .439339828 / \setdashes<2pt> \plot
7 5 5.9 2.9 5.5 1.55 / \setsolid \plot 7 5 6.1 2.3 5.5 1.55 /
\plot 7 5 6.3 1.8 6.6 1.25 / \setdashes<2pt> \plot 6.6 1.25 8.4
2.9 8.9 3.1 / \setsolid \plot 8.9 3.1 8.7 2.5 8.2 1.6 / \arrow
<6pt> [.16,.6] from  6.19 2.5 to 6.15 2.4 \arrow <6pt> [.16,.6]
from  6.31 2 to 6.29 1.9 \put {\mbox{\scriptsize
$\gamma_{\beta}$}} [lt] <0.6mm,0mm> at 6.3 2 \put
{\mbox{\scriptsize $\gamma$}} [rc] <-0.6mm,0mm> at 6.2 2.5 \put
{\mbox{\scriptsize $\delta_\beta$}} [lb] <0mm,0mm> at 8. 4 \put
{\mbox{\scriptsize
$\delta_{\alpha}\hspace{-0.03in}=\hspace{-0.02in}\gamma_{\alpha}$}}
[rc] <-1mm,0mm> at 6.2 4.2
\endpicture}
\end{center}
\caption{}\label{fig01}
\end{figure}

\vskip 5pt

Cutting $\mathcal P(\gamma)$ along $\delta_{\alpha},
\delta_{\beta}$ and $[\alpha, \beta]$ one obtains two pieces; let
the one which contains the initial part of $\gamma$ be denoted
$\mathcal P^{+}(\gamma)$. There are two simple
$\Delta_{0}$-geodesics, $\gamma_{\alpha}$ and $\gamma_{\beta}$, in
$\mathcal P(\gamma)$ such that they are asymptotic to $\alpha$ and
$\beta$ respectively, and such that their initial parts are
contained in $\mathcal P^{+}(\gamma)$. See Figure \ref{fig01}.

\vskip 10pt

\begin{lem}  Each $\Delta_{0}$-geodesic whose initial part lies in $\mathcal
P^{+}(\gamma)$ between $\gamma_{\alpha}$ and $\gamma$ or between
$\gamma$ and $\gamma_{\beta}$ is non-simple or simple-not-normal.
\end{lem}

The union of these two gaps between simple-normal
$\Delta_{0}$-geodesics formed by non-simple and simple-not-normal
$\Delta_{0}$-geodesics is called the {\bf main gap} determined by
$\gamma$.

\vskip 10pt

This lemma can be proved easily using a suitable model of the
hyperbolic plane; see \cite{zhang2004thesis} for details. The idea
is that a $\Delta_{0}$-geodesic ray whose initial part lies in
$\mathcal P^{+}(\gamma)$ between $\gamma_{\alpha}$ and $\gamma$
will not intersect $\gamma_{\alpha}$ or $\gamma$ directly, so it
must come back to intersect for first time either itself or
$\Delta_0$, hence is either non-simple or simple but not-normal
(that is, intersecting $\Delta_{0}$ obliquely). More precisely, if
$\Delta_{0}$ is a cusp or a cone point all the
$\Delta_{0}$-geodesics in the lemma are non-simple, while if
$\Delta_{0}$ is a boundary geodesic then there is a (critical)
$\Delta_{0}$-geodesic, $\rho_\gamma$, whose initial part lies in
$\mathcal P^{+}(\gamma)$ between $\gamma_{\alpha}$ and $\gamma$
such that $\rho_\gamma$ is non-simple and its only
self-intersection is at its starting point on $\Delta_{0}$ (and
hence terminates there) and it has the property that each
$\Delta_{0}$-geodesic whose initial part lies in $\mathcal
P^{+}(\gamma)$ between $\gamma_{\alpha}$ and $\rho_\gamma$ is
non-simple, while each $\Delta_{0}$-geodesic whose initial part
lies in $\mathcal P^{+}(\gamma)$ between $\rho_\gamma$ and
$\gamma$ is simple-not-normal terminating at $\Delta_{0}$. There
is a similar dichotomy for the $\Delta_{0}$-geodesics whose
initial parts lie in $\mathcal P^{+}(\gamma)$ between $\gamma$ and
$\gamma_\beta$.

\vskip 10pt

Now suppose one of $\alpha, \beta$, say $\alpha$, is a boundary
geodesic. Then there are two simple $\Delta_{0}$-geodesics in
$\mathcal P(\gamma)$ which are asymptotes to $\alpha$. They are
$\gamma_{\alpha}$ and $(-\gamma)_{\alpha}$.

\vskip 10pt

The following lemma tells us that there is an {\bf extra gap}
determined by $\gamma$ in $\mathcal P^{+}(\gamma)$ between
simple-normal $\Delta_{0}$-geodesics formed by simple-not-normal
$\Delta_{0}$-geodesics.

\vskip 10pt

\begin{lem}
Each $\Delta_{0}$-geodesic whose initial part lies in $\mathcal
P^{+}(\gamma)$ between $\delta_{\alpha}$ and $\gamma_{\alpha}$ is
simple-not-normal.
\end{lem}

This is almost self-evident from the geometry of the pair of pants
$P(\gamma)$, and is similar to the proof of the previous lemma;
see \cite{zhang2004thesis} for details.

\vskip 10pt

Note that there is a similar and symmetric picture for the
$\Delta_{0}$-geodesics whose initial parts lie in $\mathcal
P^{-}(\gamma)$.

\vskip 10pt

Hence (for non-degenerate $\gamma$) in the geometric pair of pants
$\mathcal P(\gamma)$, which is the same as $\mathcal P(-\gamma)$,
if none of $\alpha, \beta$ is a boundary geodesic then there are
two main gaps determined by $\gamma$ and $-\gamma$ respectively;
if (exactly) one of $\alpha, \beta$ is a boundary geodesic then
there are two extra gaps determined by $\gamma$ and $-\gamma$.

\vskip 10pt

The case of a degenerate simple $[\Delta_0,\Delta_0]$-geodesic
$\gamma$ is handled in a similar way. Recall that $\gamma$ is the
double cover of a $\Delta_0$-geodesic arc $\delta$ from $\Delta_0$
to an angle $\pi$ cone point $\alpha$. Then there is a simple
closed curve ${\beta}^{\prime}$, which is the boundary of a
suitable regular neighborhood of $\Delta_0 \cup \delta$ on $M$,
such that ${\beta}^{\prime}$ bounds with $\Delta_0$ and $\alpha$
an embedded (topological) pair of pants. If $\Delta_0$ is not
itself an angle $\pi$ cone point, then ${\beta}^{\prime}$ can be
realized as an interior generalized simple closed geodesic $\beta$
which bounds with $\Delta_0$ and $\alpha$ an embedded pair of
pants $\mathcal H(\Delta_0, \alpha, \beta)$ on $M$ and we can
carry out the analysis as above with suitable modifications. In
this case $\gamma$ determines no gaps if $\Delta_0$ is itself an
angle $\pi$ cone point. If $\Delta_0$ is not itself an angle $\pi$
cone point then there are two main gaps, between $\gamma$ and each
of the two $\Delta_0$-geodesics which are asymptotic to $\beta$ in
$\mathcal H(\Delta_0, \alpha, \beta)$. We say that one of the two
main gaps is determined by $\gamma$ and the other by $-\gamma$
although $\gamma=-\gamma$ in this case.

\vskip 10pt
\begin{defn} The {\bf width} of an open subinterval $\mathcal H^{\prime}$
of $\mathcal H$ is defined respectively as:\begin{itemize}
\item[(i)] $\Delta_{0}$ is a cusp: the normalized parabolic
measure, that is, the ratio of the Euclidean length of $\mathcal
H^{\prime}$ to the Euclidean length of $\mathcal H$; \item[(ii)]
$\Delta_{0}$ is a cone point: the elliptic measure, that is, the
angle (measured in radians) that $\mathcal H^{\prime}$ subtends
with respect to the cone point $\Delta_{0}$; \item[(iii)]
$\Delta_{0}$ is a boundary geodesic: the hyperbolic measure, that
is, the hyperbolic length of $\mathcal H^{\prime}$ (recall that in
this case $\mathcal H$ is the same as the distinguished boundary
geodesic $\Delta_{0}$).
\end{itemize}
\end{defn}

\vskip 10pt

\begin{defn} The {\bf combined gap} between simple-normal
$\Delta_{0}$-geodesics determined by $\gamma$ is the union of the
main gap and the extra gap (if there is any) determined by
$\gamma$. The {\bf gap function} ${\rm Gap}(\Delta_{0};\alpha,
\beta)$ when $\Delta_{0}$ is a cone point or boundary geodesic or
the {\bf normalized gap function} ${\rm
Gap}^{\prime}(\Delta_{0};\alpha, \beta)$ when $\Delta_{0}$ is a
cusp is defined as the total width of the combined gap determined
by $\gamma$, which is by symmetry the same as the total width of
the combined gap determined by $-\gamma$.
\end{defn}

We shall calculate the the gap functions in \S
\ref{s:calculation}.

\vskip 10pt

On the other hand, the following key lemma shows that the
non-simple and simple-not-normal $\Delta_{0}$-geodesics obtained
above are {\it all} the non-simple and simple-not-normal
$\Delta_{0}$-geodesics.

\vskip 10pt

\begin{lem}
Each non-simple or simple-not-normal $\Delta_{0}$-geodesic lies in
a main gap or an extra gap determined by some
$[\Delta_{0},\Delta_{0}]$-geodesic $\gamma$.
\end{lem}

\vskip 8pt

\begin{pf} First let $\delta$ be a non-simple $\Delta_{0}$-geodesic,
with its first self-intersection point $Q$, where $Q$ lies in the
geometric interior of $M$ or in $\Delta_{0}$ when $\Delta_{0}$ is
a boundary geodesic. Let $\delta_{1}$  be the part of $\delta$
from starting point to $Q$; note that $\delta_{1}$ has the shape
of a lasso. Then in the boundary of a suitable regular
neighborhood of $\delta_{1}$ there is a simple arc
$\gamma^{\prime}$ which connects $\Delta_{0}$ to itself and is
disjoint from $\delta_{1}$ (except at $\Delta_{0}$ when
$\Delta_{0}$ is a cone point); there is also a simple closed curve
$\alpha^{\prime}$ which is freely homotopic to the loop part of
$\delta_{1}$. See Figure \ref{fig02}. Let $\gamma$, $\alpha$ be
the generalized simple closed geodesics on $M$ which realize
$\gamma^{\prime}$, $\alpha^{\prime}$ in their respective free
(relative) homotopy classes in the geometric interior of $M$. An
easy geometric argument shows that $\alpha$ is disjoint from
$\delta_{1}$ and that $\gamma$ is also disjoint from $\delta_{1}$
except at $\Delta_{0}$ when $\Delta_{0}$ is a cone point or a
cusp. Furthermore, $\gamma$ and $\alpha$ cobound (together with
$\Delta_{0}$ when $\Delta_{0}$ is a boundary geodesic) an embedded
cylinder which contains $\delta_{1}$. Hence the point in $\mathcal
H$ which corresponds to the $\Delta_{0}$-geodesic $\delta$ lies in
the main gap determined by $\gamma$. See Figure \ref{fig03}

\begin{figure}
\begin{center}
\mbox{\beginpicture \setcoordinatesystem units <0.69in,0.69in>
\setplotarea x from -1.2 to 1.4 , y from -1.2 to 1.4 \circulararc
-30 degrees from 0.5 0.866 center at 0 0 \circulararc 210 degrees
from -0.5 0.866 center at 0 0 \ellipticalarc axes ratio 5:1 360
degrees from 0.5 0.866 center at 0 0.866 \ellipticalarc axes ratio
4:1 90 degrees from 0.866 0.5  center at 1.1 0.5 \ellipticalarc
axes ratio 4:1 -90 degrees from 0.866 -0.5  center at 1.1 -0.5
\ellipticalarc axes ratio 2:1 180 degrees from 0.1 -0.02  center
at 0 -0.02 \ellipticalarc axes ratio 2:1 -180 degrees from 0.15
0.02  center at 0 0.02 \plot 0.05 0.766 0.05 0.708 0. 0.705 /
\setquadratic
\plot 0 0.766 0.03 0.6 .1710100714 .461 /
\ellipticalarc axes ratio 1:1.1 -310 degrees from 0.1710100714
.461  center at 0 0
\plot -0.22 0.44 -0.05 0.51 0.085 0.53 /
\put {\mbox{\LARGE $\cdot$}} [cc] <0mm,-0.2mm> at 0.085 0.53
\ellipticalarc axes ratio 1:1.1 360 degrees from 0 0.275  center
at 0 0
\plot 0.2 0.77 0.19 0.6 0.28 0.5 / \plot -0.2 0.77 -0.19 0.6 -0.28
0.5 / \ellipticalarc axes ratio 1:1.1 -300 degrees from 0.28 0.5
center at 0 0 \put {\mbox{\LARGE $\cdots$}} [cc] <0mm,0mm> at 1.4
0 \arrow <6pt> [.16,.6] from  0.53 0 to 0.5301 -0.01 \arrow <6pt>
[.16,.6] from  0.352 0.31 to 0.364 0.292 \put {\mbox{\scriptsize
$\alpha'$}} [cb] <0mm,1mm> at 0 -0.26 \put {\mbox{\scriptsize
$\delta_{1}$}} [rc] <0mm,-1mm> at 0.35 0.3 \put {\mbox{\scriptsize
$\gamma'$}} [lc] <0mm,0mm> at 0.54 0.2 \put {\mbox{\scriptsize
$\triangle_{0}$}} [cb] <0mm,1mm> at 0 0.95
\endpicture}
\hspace{0.2in} \mbox{\beginpicture \setcoordinatesystem units
<0.69in,0.69in> \setplotarea x from -1.2 to 1.4 , y from -1.2 to
1.4 \circulararc -30 degrees from 0.5 0.866 center at 0 0
\circulararc 210 degrees from -0.5 0.866 center at 0 0
\ellipticalarc axes ratio 4:1 90 degrees from 0.866 0.5  center at
1.1 0.5 \ellipticalarc axes ratio 4:1 -90 degrees from 0.866 -0.5
center at 1.1 -0.5 \ellipticalarc axes ratio 2:1 180 degrees from
0.1 -0.02  center at 0 -0.02 \ellipticalarc axes ratio 2:1 -180
degrees from 0.15 0.02  center at 0 0.02
\setquadratic \plot 0.5 0.866 0.25 1  0 1.1 -0.25 0.99 -0.5 0.866
/
\plot 0 1.1 0.1 0.7 0.26 0.41 /
\put {\mbox{\LARGE $\cdot$}} [cc] <0mm,-0.1mm> at 0.19 0.5
\ellipticalarc axes ratio 1:1.1 270 degrees from 0 0.5  center at
0 0 \ellipticalarc axes ratio 1:1.1 -90 degrees from 0.26 0.41
center at 0 0 \ellipticalarc axes ratio 1:1.1 360 degrees from 0
0.275  center at 0 0
\plot 0 0.5 0.11 0.507 0.2 0.5 /
\plot 0 1.1 0.15 0.67 0.28 0.5 / \plot 0 1.1 -0.15 0.67 -0.28 0.5
/ \ellipticalarc axes ratio 1:1.1 -300 degrees from 0.28 0.5
center at 0 0 \put {\mbox{\LARGE $\cdots$}} [cc] <0mm,0mm> at 1.4
0 \arrow <6pt> [.16,.6] from  0.53 0 to 0.5301 -0.01 \arrow <6pt>
[.16,.6] from  0.352 0.31 to 0.364 0.292
\put {\mbox{\scriptsize $\alpha'$}} [cb] <0mm,1mm> at 0 -0.26 \put
{\mbox{\scriptsize $\delta_{1}$}} [rc] <0mm,-1mm> at 0.35 0.3 \put
{\mbox{\scriptsize $\gamma'$}} [lc] <0mm,0mm> at 0.54 0.2
\put {\mbox{\scriptsize $\triangle_{0}$}} [cb] <0mm,1mm> at 0 1.1
\endpicture}
\end{center}
\caption{}\label{fig02}
\end{figure}

\begin{figure}
\begin{center}
\mbox{\beginpicture \setcoordinatesystem units <0.24in,0.24in>
\setplotarea x from 0 to 9, y from 0 to 6 \ellipticalarc axes
ratio 2:1 360 degrees from 3.6 5 center at 5 5 \setquadratic
\setdashes<2pt> \plot 9.560660172  2.560660172  9.399493762
2.618001774  9.150278635  2.565905573  8.837409705  2.409471107
8.491512762  2.164011274  8.146446610  1.853553390  7.835988726
1.508487238  7.590528893  1.162590295  7.434094427  .8497213651
7.381998226  .6005062381  7.439339828  .439339828 / \plot
2.560660172  .439339828  2.618001774  .6005062381  2.565905573
.8497213651  2.409471107  1.162590295  2.164011274  1.508487238
1.853553390  1.853553390  1.508487238  2.164011274  1.162590295
2.409471107  .8497213651  2.565905573  .6005062381  2.618001774
.439339828  2.560660172 / \setsolid \plot 9.560660172  2.560660172
9.618001774  2.399493762  9.565905573  2.150278635  9.409471107
1.837409705  9.164011274  1.491512762  8.853553390  1.146446610
8.508487238  .8359887260  8.162590295  .5905288924  7.849721365
.4340944273  7.600506238  .3819982259  7.439339828  .439339828 /
\plot 2.560660172  .439339828  2.399493762  .3819982259
2.150278635  .4340944273  1.837409705  .5905288924  1.491512762
.8359887260  1.146446610  1.146446610  .8359887260  1.491512762
.5905288924  1.837409705  .4340944273  2.150278635  .3819982259
2.399493762  .439339828  2.560660172 /
\plot 9.560660172  2.560660172 7.5 4 6.4 5 / \plot .439339828
2.560660172 2.5 4 3.6 5 / \plot 2.560660172  .439339828 5 1
7.439339828  .439339828 /
\plot 5.224157199 4.3 5.1 1.7 5 1  / \setdashes<2pt> \plot 5 1 4.8
2  4.775842801 5.7 /
\setsolid \plot 4.628894061 5.66 4.62 5.5 4.76 5.5 / \plot 5.39
4.3 5.39 4.15 5.21 4.14 /
\plot 4.25 4.4 4.1 1.7 3.75 0.85 / \setdashes<2pt> \plot 3.75 0.85
2.4 2.1 2 3.6 / \setsolid \plot 2 3.6 3 2.6  3.5 0.78 /
\setdashes<2pt>
\plot 3.5 0.78 2.3 1.5 1.2 3 / \setsolid \plot 1.2 3 1.6 3.1 3.2
2.6 / \plot 4.4 4.35 4.4 4.2 4.25 4.23 / \put {\mbox{\LARGE
$\cdot$}} [cc] <0mm,0mm> at 2.93 2.66 \arrow <6pt> [.16,.6] from
5.19 3.7 to 5.188 3.6 \arrow <6pt> [.16,.6] from  4.23 3.7 to 4.22
3.6 \put {\mbox{\small $\triangle_0$}} [cb] <0mm,2mm> at 5 5.5
\put {\mbox{\small $\alpha$}} [rt] <-1mm,-1.6mm> at 1.4 1.2
\put {\mbox{\small $\gamma$}} [lc] <1mm,0mm> at 5.1 3.2 \put
{\mbox{\small $\delta$}} [rc] <-1mm,0mm> at 4.3 3.2
\endpicture}
\hspace{0.4in} \mbox{\beginpicture \setcoordinatesystem units
<0.29in,0.29in> \setplotarea x from 3.5 to 9, y from 0.5 to 6
\plot 9.560660172 2.560660172 7 5 4 2  / \setquadratic \plot
9.560660172 2.560660172 9.639852575  2.377642961  9.607468267
2.108715941 9.466677248 1.780203564  9.231261125  1.424262911
8.924264069 1.075735931 8.575737089  .7687388748  8.219796436
.5333227521 7.891284059 .3925317336  7.622357039  .3601474247
7.439339828 .439339828 / \setdashes<2pt> \plot 9.560660172
2.560660172 9.377642961 2.639852575  9.108715941  2.607468267
8.780203564 2.466677248 8.424262911  2.231261125  8.075735931
1.924264069 7.768738875 1.575737089  7.533322752  1.219796435
7.392531733 .8912840588 7.360147425  .6223570393  7.439339828
.439339828 / \put {\mbox{\LARGE $\cdot$}} [cc] <0mm,0mm> at 7 5
\put {\mbox{\scriptsize $\triangle_{0}$}} [cb] <0mm,1mm> at 7 5
\put {\mbox{\LARGE $\cdot$}} [cc] <0mm,0mm> at 4 2 \put
{\mbox{\scriptsize $\alpha$}} [rt] <-1mm,0mm> at 4 2 \setsolid
\plot 4 2 6.3 1.2 7.439339828  .439339828 / \setdashes<2pt> \plot
5.5 1.55 5.1 2 5 3 /
\setsolid \plot 5 3 5.3 2.3 5.2 1.62 / \setdashes<2pt> \plot 5.2
1.62 4.8 2  4.7 2.7 / \setsolid \plot 4.7 2.7 5 2.7 5.5 2.3 / \put
{\mbox{\LARGE $\cdot$}} [cc] <0mm,0mm> at 5.25 2.5
\setsolid \plot 7 5 6.1 2.3 5.5 1.55 /
\plot 6.6 1.25 7.1 2.5 7 5 / \setdashes<2pt> \plot 7 5 6.3 1.8 6.6
1.25 /
\setsolid \arrow <6pt> [.16,.6] from  6.19 2.5 to 6.15 2.4
\arrow <6pt> [.16,.6] from  7.02 2 to 7 1.9 \put
{\mbox{\scriptsize $\gamma$}} [lt] <0.6mm,0mm> at 7 2 \put
{\mbox{\scriptsize $\delta$}} [rc] <-0.6mm,0mm> at 6.2 2.5
\endpicture}
\end{center}
\caption{}\label{fig03}
\end{figure}

\vskip 6pt

Next let $\delta$ be a simple-not-normal $\Delta_{0}$-geodesic
which terminates at $\Delta_{0}$ itself; in this case $\Delta_{0}$
is a boundary geodesic and $\mathcal H$ is $\Delta_{0}$ itself.
Then the boundary of a suitably chosen regular neighborhood of
$\delta \cup \mathcal H$ consist of two disjoint simple closed
curves in the geometric interior of $M$. Let their geodesic
realizations be (disjoint) generalized simple closed geodesics
$\alpha$ and $\beta$. Then $\alpha, \beta$ bound with $\Delta_{0}$
an embedded pair of pants which contains $\delta$ in a main gap
determined by the $[\Delta_{0}, \Delta_{0}]$-geodesic $\gamma$
which is the geodesic realization of $\delta$ in its free relative
homotopy class.

\vskip 6pt

Finally let $\delta$ be a simple-not-normal $\Delta_{0}$-geodesic
which terminates at a boundary geodesic $\Delta_{1}$ which is
different from $\Delta_{0}$. The boundary of suitably chosen
regular neighborhood of $\delta \cup \Delta_{1}$ on $M$ is a
simple arc connecting $\Delta_{0}$ to itself and is disjoint from
$\delta$. Its geodesic realization is a $[\Delta_{0},
\Delta_{0}]$-geodesic, $\gamma$, which is disjoint from $\delta$.
Now $\Delta_{1}, \gamma$ bound with $\Delta_{0}$ an embedded
cylinder which contains $\delta$. Hence $\delta$ lies in the extra
gap determined by $\gamma$ or $-\gamma$.
\end{pf}

\vskip 10pt

\vskip 30pt
\section{{\bf Calculating the gap functions}}\label{s:calculation}
\vskip 30pt

In this section we calculate the gap function ${\rm
Gap}(\Delta_{0};\alpha, \beta)$ when $\Delta_{0}$ is a cone point
or a boundary geodesic, it is the width of the combined gap
determined by a simple $[\Delta_{0},\Delta_{0}]$-geodesic $\gamma$
on $M$.

\vskip 10pt

Recall that $\alpha, \beta$ are the generalized simple closed
geodesics determined by $\gamma$ and $\mathcal P(\gamma)$ is the
geometric pair of pants that $\alpha, \beta$ bound with
$\Delta_{0}$ on $M$.

\vskip 20pt

{\it Case} 1.  $\Delta_{0}$ is a cone point of cone angle $\theta
\in (0, \pi]$.

\vskip 10pt

In this case the width of the main gap determined by $\gamma$ is
the angle between $\gamma_{\alpha}$ and $ \gamma_{\beta}$.

\vskip 10pt

Let $x$ be the angle between $\delta_{\alpha}$ and
$\gamma_{\alpha}$ and let $y$ be the angle between
$\delta_{\beta}$ and $\gamma_{\beta}$.

\vskip 20pt

{\it Subcase} 1.1. Both $\alpha$ and $\beta$ are interior
generalized simple closed curves.

\vskip 10pt

In this case the width of the combined gap determined by $\gamma$
is the angle between $\gamma_{\alpha}$ and $ \gamma_{\beta}$ and
is equal to $\frac{\theta}{2} - (x+y)$.

\vskip 10pt

By a formula in Fenchel \cite{fenchel1989book} VI.3.2 (line 10,
page 87),
\begin{eqnarray}
\sinh |\delta_{\alpha}| =
\frac{\cosh\frac{|\beta|}{2}+\cos\frac{\theta}{2}\cosh\frac{|\alpha|}{2}}
{\sin\frac{\theta}{2}\sinh\frac{|\alpha|}{2}},
\end{eqnarray}
\begin{eqnarray}
\sinh |\delta_{\beta}| =
\frac{\cosh\frac{|\alpha|}{2}+\cos\frac{\theta}{2}\cosh\frac{|\beta|}{2}}
{\sin\frac{\theta}{2}\sinh\frac{|\beta|}{2}}.
\end{eqnarray}
Hence
\begin{eqnarray}\label{eqn:1.1.x}
\tan x = \frac{1}{\sinh |\delta_{\alpha}|} =
\frac{\sin\frac{\theta}{2}\sinh\frac{|\alpha|}{2}}
{\cosh\frac{|\beta|}{2}+\cos\frac{\theta}{2}\cosh\frac{|\alpha|}{2}},
\end{eqnarray}
\begin{eqnarray}\label{eqn:1.1.y}
\tan y = \frac{1}{\sinh |\delta_{\beta}|} =
\frac{\sin\frac{\theta}{2}\sinh\frac{|\beta|}{2}}
{\cosh\frac{|\alpha|}{2}+\cos\frac{\theta}{2}\cosh\frac{|\beta|}{2}}.
\end{eqnarray}
From these one can derive that
\begin{eqnarray}
\tan \, (x+y) =
\frac{\sin\frac{\theta}{2}\sinh\frac{|\alpha|+|\beta|}{2}}
{1+\cos\frac{\theta}{2}\cosh\frac{|\alpha|+|\beta|}{2}}
\end{eqnarray}
and hence that
\begin{eqnarray}
\tan \frac {x+y}{2} = \tan \frac {\theta}{4} \tanh \frac
{|\alpha|+|\beta|}{4}.
\end{eqnarray}
Thus
\begin{eqnarray}
\tan \left(\frac {\theta}{4} - \frac {x+y}{2}\right) &=& \frac
{\tan \frac{\theta}{4} \left(1-\tanh
\frac{|\alpha|+|\beta|}{4}\right)} {1+ \tan^{2}\frac{\theta}{4} \,
\tanh \frac{|\alpha|+|\beta|}{4}}\\
&=& \frac{\sin\frac{\theta}{2}}{\cos\frac{\theta}{2}+
\exp\frac{|\alpha|+|\beta|}{2}}.
\end{eqnarray}
Hence in this case we have
\begin{eqnarray*}
{\rm Gap}(\Delta_{0};\alpha, \beta) &=& \frac{\theta}{2} - (x+y)\\
&=& 2 \tan^{-1}
\left(\frac{\sin\frac{\theta}{2}}{\cos\frac{\theta}{2}+
\exp\frac{|\alpha|+|\beta|}{2}} \right ).
\end{eqnarray*}

\vskip 10pt

\begin{figure}
\begin{center}
\mbox{\beginpicture \setcoordinatesystem units <1in,1in>
\setplotarea x from -1.2 to 1.2, y from -1.1 to 1.2 \circulararc
360 degrees from 1 0 center at 0 0 \startrotation by -.7071067810
-.7071067810 about 0 0 \plot .5714285714  0 0 0 0 .5714285714  /
\setdashes<2pt> \circulararc -59.48976262 degrees from .5714285714
0 center at 1.160714286 0 \circulararc 59.48976262 degrees from 0
.5714285714  center at 0 1.160714286 \circulararc 55.15983610
degrees from .3693509756  .3693509756 center at  .8615384618
.8615384618 \circulararc -55.15983610 degrees from .3693509756
.3693509756 center at  .8615384618 .8615384618 \setsolid
\circulararc -20.48 degrees from .5714285714 0 center at
1.160714286 0 \circulararc 20.48 degrees from 0 .5714285714
center at 0 1.160714286 \circulararc 25.15983610 degrees from
.3693509756  .3693509756 center at  .8615384618 .8615384618
\circulararc -25.15983610 degrees from .3693509756  .3693509756
center at  .8615384618 .8615384618
\plot .52 0 .52 .05 .573 .05 / \plot 0 .52 .05 .52 .05 .573 /
\plot 0.598 .1749579531 0.55 .192 0.57 .2294764655 / \plot
.1749579531 0.598 .192 0.55 .2294764655 0.57 /
\plot 0 0 .48 .2828571428 / \plot 0 0 .2828571428 .49 /
\setdashes<2pt> \plot .48 .2828571428 .8615384616 .5076923076 /
\plot .2828571428 .49 .5076923076 .8615384616 / \put
{\mbox{\scriptsize $\triangle_{0}$}} [cb] <0mm,0mm> at 0 0 \put
{\mbox{\scriptsize $\alpha$}} [rc] <-3mm,-1mm> at 0.5 0.18 \put
{\mbox{\scriptsize $\beta$}} [lc] <3mm,-1mm> at 0.18 0.5 \put
{\mbox{\scriptsize $\delta_{\alpha}$}} [rb] <-5mm,0mm> at 0.2
.1178571428 \put {\mbox{\scriptsize $\delta_{\beta}$}} [lb]
<5mm,0mm> at .1178571428 0.2 \put {\mbox{\scriptsize
$\gamma_{\alpha}$}} [rt] <-0.4mm,0mm> at 0.3 .1767857142 \put
{\mbox{\scriptsize $\gamma_{\beta}$}} [lt] <0.4mm,0mm> at
.1767857142 0.3
\stoprotation
\put {\mbox{\small \sc Figure 4. \rm Subcase
1.1\label{fig04.1.1}}} [ct] <0mm,0mm> at 0 -1.3
\endpicture}
\hspace{0.2in} \raisebox{-5ex}{\mbox{\beginpicture
\setcoordinatesystem units <0.2in,0.2in> \setplotarea x from -1 to
10, y from 1 to 6 \plot 0 0  5 5 9 1 6 -2 / \plot 5 5 4.8 -1.2 /
\plot 8.8 1.2 8.6 1 8.8 0.8 / \plot 6.2 -1.8 6 -1.6 5.82 -1.8 /
\setquadratic \plot 0 0 3.2 -0.5 6 -2 / \put {\mbox{\LARGE
$\cdot$}} [cc] <0mm,0mm> at 5 5 \put {\mbox{\scriptsize
$\triangle_{0}$}} [cb] <0mm,1mm> at 5 5 \put {\mbox{\LARGE
$\cdot$}} [cc] <0mm,0mm> at 0 0 \put {\mbox{\scriptsize $\alpha$}}
[rt] <-1mm,0mm> at 0 0 \put {\mbox{\scriptsize $\gamma_{\beta}$}}
[lc] <0.6mm,0mm> at 5 2 \put {\mbox{\scriptsize $\beta$}} [lt]
<2mm,0mm> at 7.2 -0.5 \put {\mbox{\scriptsize $\delta_\beta$}}
[lb] <2mm,0mm> at 6.4 3.3 \put {\mbox{\scriptsize
$\delta_{\alpha}\hspace{-0.03in}=\hspace{-0.03in}\gamma_{\alpha}$}}
[rb] <-1mm,0mm> at 2.5 2.5 \put {\mbox{\small \sc Figure 5. \rm
Subcase 1.3 \label{fig05.1.3}}} [ct] <0mm,0mm> at 5 -5
\endpicture}}
\end{center}
\end{figure}

{\it Subcase} 1.2. $\alpha$ is a boundary geodesic and $\beta$ is
an interior generalized simple closed geodesic.

\vskip 10pt

In this case the width of the combined gap determined by $\gamma$
is the angle between $\delta_{\alpha}$ and $ \gamma_{\beta}$ and
is equal to $\frac{\theta}{2} - y$. Hence by (\ref{eqn:1.1.y}) we
have
\begin{eqnarray} {\rm Gap}(\Delta_{0};\alpha, \beta)=
\frac{\theta}{2} - \tan^{-1} \left (
\frac{\sin\frac{\theta}{2}\sinh\frac{|\beta|}{2}
}{\cosh\frac{|\alpha|}{2}+\cos\frac{\theta}{2}\cosh\frac{|\beta|}{2}
} \right ).
\end{eqnarray}

\vskip 10pt

{\it Subcase} 1.3. $\alpha$ is a cone point of cone angle $\varphi
\in (0, \pi]$ and $\beta$ is an interior generalized simple closed
geodesic.

\vskip 10pt

Note that in this case $\gamma_{\alpha}$ coincides with
$\delta_{\alpha}$ and hence $x=0$. Therefore the width of the
combined gap determined by $\gamma$ is the angle between
$\delta_{\alpha}$ and $ \gamma_{\beta}$ and is equal to
$\frac{\theta}{2} - y$.

\vskip 10pt

Now by a formula in Fenchel \cite{fenchel1989book} VI.3.3 (line
13, page 88),
\begin{eqnarray}\label{}
\sinh |\delta_{\beta}| =
\frac{\cos\frac{\varphi}{2}+\cos\frac{\theta}{2}\cosh\frac{|\beta|}{2}}
{\sin\frac{\theta}{2}\sinh\frac{|\beta|}{2}}.
\end{eqnarray}
Hence
\begin{eqnarray}\label{eqn:1.3.y}
\tan y = \frac{1}{\sinh |\delta_{\beta}|} =
\frac{\sin\frac{\theta}{2}\sinh\frac{|\beta|}{2}}
{\cos\frac{\varphi}{2}+\cos\frac{\theta}{2}\cosh\frac{|\beta|}{2}}.
\end{eqnarray}
Thus in this case we have
\begin{eqnarray} {\rm
Gap}(\Delta_{0};\alpha, \beta)= \frac{\theta}{2} - \tan^{-1} \left
( \frac{\sin\frac{\theta}{2}\sinh\frac{|\beta|}{2}
}{\cos\frac{\varphi}{2}+\cos\frac{\theta}{2}\cosh\frac{|\beta|}{2}}
\right ).
\end{eqnarray}

\vskip 20pt

{\it Case} 2.  $\Delta_{0}$ is a boundary geodesic of length
$l>0$.

\vskip 10pt

In this case the width of the main gap determined by $\gamma$ is
the distance between $\gamma_{\alpha}$ and $ \gamma_{\beta}$ along
$\Delta_{0}$.

\vskip 10pt

Let $x$ be the distance between $\delta_{\alpha}$ and
$\gamma_{\alpha}$ along $\Delta_{0}$ and let $y$ be the distance
between $\delta_{\beta}$ and $\gamma_{\beta}$ along $\Delta_{0}$.

\vskip 10pt

We shall see that all calculations in this case are parallel to
those in Case 1.

\vskip 20pt

{\it Subcase} 2.1. Both $\alpha$ and $\beta$ are interior
generalized simple closed curves.

\vskip 10pt

In this case the width of the combined gap determined by $\gamma$
is the distance between $\gamma_{\alpha}$ and $ \gamma_{\beta}$
along $\Delta_{0}$ and is equal to $\frac{l}{2} - (x+y)$.

\vskip 10pt

By the cosine rule for right angled hexagons on the hyperbolic
plane (c.f. Fenchel \cite{fenchel1989book} VI.3.1, page 86, or
Beardon \cite{beardon1983book} Theorem 7.19.2, page 161),
\begin{eqnarray}
\cosh |\delta_{\alpha}| =
\frac{\cosh\frac{|\beta|}{2}+\cosh\frac{l}{2}\cosh\frac{|\alpha|}{2}}
{\cosh\frac{l}{2}\sinh\frac{|\alpha|}{2}},
\end{eqnarray}
\begin{eqnarray}
\cosh |\delta_{\beta}| =
\frac{\cosh\frac{|\alpha|}{2}+\cosh\frac{l}{2}\cosh\frac{|\beta|}{2}}
{\sinh\frac{l}{2}\sinh\frac{|\beta|}{2}}.
\end{eqnarray}
Hence
\begin{eqnarray}\label{eqn:2.1.x}
\tanh x = \frac{1}{\cosh |\delta_{\alpha}|} =
\frac{\sinh\frac{l}{2}\sinh\frac{|\alpha|}{2}}
{\cosh\frac{|\beta|}{2}+\cosh\frac{l}{2}\cosh\frac{|\alpha|}{2}},
\end{eqnarray}
\begin{eqnarray}\label{eqn:2.1.y}
\tanh y = \frac{1}{\cosh |\delta_{\beta}|} =
\frac{\sinh\frac{l}{2}\sinh\frac{|\beta|}{2}}
{\cosh\frac{|\alpha|}{2}+\cosh\frac{l}{2}\cosh\frac{|\beta|}{2}}.
\end{eqnarray}
From these one can derive that
\begin{eqnarray}
\tanh \, (x+y) =
\frac{\sinh\frac{l}{2}\sinh\frac{|\alpha|+|\beta|}{2}}
{1+\cosh\frac{l}{2}\cosh\frac{|\alpha|+|\beta|}{2}}
\end{eqnarray}
and hence that
\begin{eqnarray}
\tanh \frac {x+y}{2} = \tanh \frac {l}{4} \tanh \frac
{|\alpha|+|\beta|}{4}.
\end{eqnarray}
Thus
\begin{eqnarray}
\tanh \left(\frac {l}{4} - \frac {x+y}{2}\right) &=& \frac {\tanh
\frac{l}{4} \left(1-\tanh \frac{|\alpha|+|\beta|}{4}\right)} {1-
\tanh^{2}\frac{l}{4} \, \tanh \frac{|\alpha|+|\beta|}{4}}\\
&=&\frac {\sinh \frac{l}{2}}{\cosh \frac{l}{2}+
\exp\frac{|\alpha|+|\beta|}{2}}.
\end{eqnarray}
Hence in this case we have
\begin{eqnarray*}
{\rm Gap}(\Delta_{0};\alpha, \beta) &=& \frac{l}{2} - (x+y)\\
&=& 2\tanh^{-1} \left ( \frac {\sinh \frac{l}{2}}
{\cosh\frac{l}{2}+ \exp\frac{|\alpha|+|\beta|}{2}}\right ).
\end{eqnarray*}

\vskip 10pt

\begin{figure}
\begin{center}
\mbox{\beginpicture \setcoordinatesystem units <1in,1in>
\setplotarea x from -1.2 to 1.2, y from -1.3 to 1.2 \circulararc
360 degrees from 1 0 center at 0 0 \startrotation by -1 0 about 0
0 \plot  .5714285714  0 0 0 / \setdashes<2pt> \circulararc
-59.48976262 degrees from .5714285714 0 center at 1.160714286 0
\circulararc 59.48976262 degrees from 0 .5714285714  center at 0
1.160714286 \circulararc 55.15983610 degrees from .3693509756
.3693509756 center at  .8615384618 .8615384618 \circulararc
-55.15983610 degrees from .3693509756  .3693509756 center at
.8615384618 .8615384618 \setsolid \circulararc -20.48 degrees from
.5714285714 0 center at 1.160714286 0 \circulararc 20.48 degrees
from 0 .5714285714  center at 0 1.160714286 \circulararc
25.15983610 degrees from .3693509756  .3693509756 center at
.8615384618 .8615384618 \circulararc -25.15983610 degrees from
.3693509756  .3693509756 center at  .8615384618 .8615384618
\plot .52 0 .52 .05 .573 .05 /
\plot 0.598 .1749579531 0.55 .192 0.57 .2294764655 / \plot
.1749579531 0.598 .192 0.55 .2294764655 0.57 /
\circulararc -5.8 degrees from .088894345 0 center at 5.669101893
0 \plot 0.14 0 0.14 0.05 0.1 0.05 / \setdashes<2pt> \circulararc
4.2 degrees from .1763947833 .9843195012 center at 5.669101893 0
\setsolid
\put {\mbox{\scriptsize $\triangle_{0}$}} [cb] <0mm,1mm> at 0 0

\put {\mbox{\scriptsize $\alpha$}} [lc] <3mm,-1mm> at 0.52 0.43
\put {\mbox{\scriptsize $\delta_{\alpha}$}} [rb] <-5mm,0mm> at 0.4
.1178571428
\put {\mbox{\scriptsize $\gamma_{\alpha}$}} [rt] <-0.4mm,0mm> at
0.09 .1767857142
\stoprotation
\startrotation by 0 -1 about 0 0 \plot 0 0 0 .5714285714 /
\setdashes<2pt> \circulararc -59.48976262 degrees from .5714285714
0 center at 1.160714286 0 \circulararc 59.48976262 degrees from 0
.5714285714  center at 0 1.160714286 \circulararc 55.15983610
degrees from .3693509756  .3693509756 center at  .8615384618
.8615384618 \circulararc -55.15983610 degrees from .3693509756
.3693509756 center at  .8615384618 .8615384618 \setsolid
\circulararc -20.48 degrees from .5714285714 0 center at
1.160714286 0 \circulararc 20.48 degrees from 0 .5714285714
center at 0 1.160714286 \circulararc 25.15983610 degrees from
.3693509756  .3693509756 center at  .8615384618 .8615384618
\circulararc -25.15983610 degrees from .3693509756  .3693509756
center at  .8615384618 .8615384618 \put {\mbox{\scriptsize
$\beta$}} [lb] <5mm,-3mm> at 0.4 .1178571428 \put
{\mbox{\scriptsize $\delta_{\beta}$}} [lb] <0mm,0mm> at .16 0.6
\put {\mbox{\scriptsize $\gamma_{\beta}$}} [lt] <0.4mm,0mm> at 0.2
0.1
\plot 0 .52 .05 .52 .05 .573 / \plot 0.598 .1749579531 0.55 .192
0.57 .2294764655 / \plot .1749579531 0.598 .192 0.55 .2294764655
0.57 / \stoprotation \startrotation by -1 0 about 0 0 \circulararc
5.8 degrees from -.088894345 0 center at -5.669101893 0 \plot
-0.14 0 -0.14 0.05 -0.1 0.05 / \setdashes<2pt> \circulararc -4.2
degrees from -.1763947833 .9843195012 center at -5.669101893 0
\stoprotation \put {\mbox{\small \sc Figure 6. \rm Subcase 2.1\ \
\label{fig06}}} [ct] <0mm,0mm> at 0 -1.2
\endpicture}
\hspace{0.2in} \raisebox{-13ex}{\mbox{\beginpicture
\setcoordinatesystem units <0.05in,0.05in> \setplotarea x from 0
to 40, y from 0 to 26 \plot 2 9 12 27 / \plot 28 27 35 10 25 5 /
\plot 20 25 20 9 / \plot 13 26.5 12.5 25.6 11.5 26.1 / \plot 21 25
21 24 20 24 / \plot 27 26.5 27.4 25.4 28.5 25.8 / \plot 34 9.5
33.5 10.6 34.5 11 / \plot 26 5.5 25. 6.5 24 5.9 / \setquadratic
\plot 2 9 14 11 25 5 / \plot 12 27 20 25 28 27 / \put
{\mbox{\LARGE $\cdot$}} [cc] <0mm,0mm> at 2 9 \put
{\mbox{\scriptsize $\alpha$}} [rt] <-1mm,0mm> at 2 9 \put
{\mbox{\scriptsize $\delta_{\beta}$}} [lb] <0mm,0mm> at 32.5 18
\put {\mbox{\scriptsize $\gamma_{\beta}$}} [lb] <0mm,0mm> at 20.5
17 \put {\mbox{\scriptsize
$\delta_{\alpha}\hspace{-0.03in}=\hspace{-0.02in}\gamma_{\alpha}$}}
[rc] <-1mm,0mm> at 7 18 \put {\mbox{\scriptsize $\triangle_{0}$}}
[cb] <0mm,1mm> at 20 25 \put {\mbox{\small \sc Figure 7. \rm
Subcase 2.3 \label{fig07}}} [ct] <0mm,0mm> at 18 -7.2
\endpicture}}
\end{center}
\end{figure}

{\it Subcase} 2.2. $\alpha$ is a boundary geodesic and $\beta$ is
an interior generalized simple closed geodesic.

\vskip 10pt

In this case the width of the combined gap determined by $\gamma$
is the distance between $\delta_{\alpha}$ and $ \gamma_{\beta}$
along $\Delta_{0}$ and is equal to $\frac{l}{2} - y$. Hence by
(\ref{eqn:2.1.y}) we have
\begin{eqnarray} {\rm
Gap}(\Delta_{0};\alpha, \beta)= \frac{l}{2} - \tanh^{-1} \left (
\frac{\sinh\frac{l}{2}\sinh\frac{|\beta|}{2}
}{\cosh\frac{|\alpha|}{2}+\cosh\frac{l}{2}\cosh\frac{|\beta|}{2} }
\right ).
\end{eqnarray}

\vskip 10pt

{\it Subcase} 2.3. $\alpha$ is a cone point of cone angle $\varphi
\in (0, \pi]$ and $\beta$ is an interior generalized simple closed
geodesic.

\vskip 10pt

Note that in this case $\gamma_{\alpha}$ coincides with
$\delta_{\alpha}$ and hence $x=0$. Hence the width of the combined
gap determined by $\gamma$ is the distance between
$\delta_{\alpha}$ and $ \gamma_{\beta}$ along $\Delta_{0}$ and is
equal to $\frac{l}{2} - y$.

\vskip 10pt

Now by a formula in Fenchel \cite{fenchel1989book} VI.3.2 (line 8,
page 87),
\begin{eqnarray}
\cosh |\delta_{\beta}| =
\frac{\cos\frac{\varphi}{2}+\cosh\frac{l}{2}\cosh\frac{|\beta|}{2}}
{\sinh\frac{l}{2}\sinh\frac{|\beta|}{2}}.
\end{eqnarray}
Hence
\begin{eqnarray}\label{eqn:2.3.y}
\tanh y = \frac{1}{\cosh |\delta_{\beta}|} =
\frac{\sinh\frac{l}{2}\sinh\frac{|\beta|}{2}}
{\cos\frac{\varphi}{2}+\cosh\frac{l}{2}\cosh\frac{|\beta|}{2}}.
\end{eqnarray}
Thus in this case we have
\begin{eqnarray} {\rm
Gap}(\Delta_{0};\alpha, \beta)= \frac{l}{2} - \tanh^{-1} \left (
\frac{\sinh\frac{l}{2}\sinh\frac{|\beta|}{2}
}{\cos\frac{\varphi}{2}+\cosh\frac{l}{2}\cosh\frac{|\beta|}{2}}
\right ).
\end{eqnarray}

\vskip 15pt

\begin{rmk}
We remark that the formulas in Case 0 for the normalized width
${\rm Gap}^{\prime}(\Delta_{0};\alpha, \beta)$ when $\Delta_{0}$
is a cusp can be derived by similar (and simpler) calculations or
by considering the first order infinitesimal terms of those
formulas with respect to $\theta$ in Case 1 or with respect to $l$
in Case 2. Hence all derivations in Case 0 are omitted.
\end{rmk}

\vskip 30pt
\section{{\bf Generalization of the Birman--Series Theorem}}\label{s:gBS}
\vskip 30pt

The celebrated Birman--Series Theorem \cite{birman-series1985t} in
its simplest form states that complete simple geodesics on a
closed hyperbolic surface are sparsely distributed.

\vskip 10pt

More precisely, let $M$ be a hyperbolic surface possibly with
boundary such that $M$ is either compact or obtained from a
compact surface by removing a finite set of points which form the
cusps and such that each boundary component of $M$ is a simple
closed geodesic. A geodesic on $M$ is said to be {\it complete} if
it is either closed and smooth, or open and of infinite length in
both directions. Hence a complete geodesic never intersects
$\partial M$. Let $G_{k}$ be the family of complete geodesics on
$M$ which have at most $k$, counted with multiplicity, transversal
self-intersections, $k \ge 0$. Then the main result in
\cite{birman-series1985t} is:

\vskip 8pt

\begin{thm} For each $k \ge 0$, the point set $S_{k}$ which is the
union of all geodesics, as point sets, in $G_{k}$ is nowhere dense
and has Hausdorff dimension one.
\end{thm}

\vskip 10pt

In this section we show that this theorem extends to the case when
$M$ is a compact hyperbolic cone-surface with geometric boundary
where each cone point has cone angle in $(0, \pi]$, with complete
geodesics replaced by complete-normal ones. This is the set of
geodesics which are either complete, or intersect the boundary
perpendicularly.

\vskip 8pt

\begin{thm} Let $M$ be a compact hyperbolic cone-surface with
geometric boundary where each cone point has cone angle in $(0,
\pi]$, and let $G_{k}$ be the family of complete-normal geodesics
on $M$ which have at most $k$ transversal self-intersections, $k
\ge 0$. Then, for each $k \ge 0$, the point set $S_{k}$ which is
the union of all geodesics, as point sets, in $G_{k}$ is nowhere
dense and has Hausdorff dimension one.
\end{thm}

\vskip 10pt

The proof of this generalization is is essentially the same as
that of the original Birman--Series theorem given in
\cite{birman-series1985t}. Hence for simplicity we shall only
sketch the proof of the theorem for the case $k=0$, that is, for
simple complete-normal geodesics; the reader is referred to
\cite{birman-series1985t} for omitted details.

\vskip 10pt

We only need to consider the case where $M$ has no geodesic
boundary components; for if $M$ has nonempty geodesic boundary we
can replace $M$ by the double of $M$ along its geodesic boundary.
We also assume for clarity that each cone point of $M$ has cone
angle less than $\pi$. We decompose the set $G_0$ into finitely
many subsets and prove the conclusion for each such subset. For
the subset of simple complete geodesics on $M$, that is, the
geodesics which never start from or terminate at cusps or cone
points, the proof is the same as that in \cite{birman-series1985t}
with little modification (which can be seen from the sketch
below). For the subset of simple normal geodesics which connect a
given cusp or cone point to another (possibly the same) given cusp
or cone point, it is easy to see that in this subset each such
geodesic is isolated in suitable neighborhoods of its endpoints
and hence the conclusion follows. Thus it remains to prove the
conclusion for the subset of simple complete-normal geodesics
which starts from a given cusp or cone point $P$ and never
terminates at any geometric boundary component.

\vskip 10pt

One can cut $M$ along normal geodesics connecting cusps or cone
points to form a (convex) fundamental polygon $R$ for $M$ in the
hyperbolic plane. Let $A=\{ a_1, a_2, \cdots, a_m \}$ denote the
ordered set of vertices and oriented sides of $R$ with
anti-clockwise ordering with some arbitrary but henceforth fixed
initial element $a_1$.

\vskip 10pt

Let $J_0$ be the set of oriented simple-normal geodesic arcs
$\gamma$ on $M$ such that the initial point and the ending point
of $\gamma$ lie in $\partial R$. (Note that except at its initial
point or ending point $\gamma$ cannot pass through a vertex of
$R$.) For $\gamma \in J_0$, we call the components of $\gamma \cap
R$  the {\it segments} of $\gamma$ and the points of $\gamma \cap
\partial R$  the partition points of $\gamma$. We label the
partition points $t_0, t_1, \cdots, t_n$ in the order in which
they occur along $\gamma$ (note that we treat $t_i \in \partial R$
as the initial point of the segment of $\gamma$ from $t_i$ to
$t_{i+1}$) and we set $\parallel \gamma
\parallel = n$ as the combinatorial length of $\gamma$.

For $\gamma \in J_0$, the segments of $\gamma$ give rise to a {\it
simple diagram} on $R$ which is a collection of finitely many
pairwise disjoint (geodesic) arcs joining pairs of distinct
elements of $A$. Two simple diagrams are regarded as being
identical if they agree up to isotopy supported on each side of
$R$. For $a_i, a_j \in A, i \neq j$, let $n_{ij}$ denote the
number of arcs joining $a_i$ to $a_j$ in the given simple diagram.
The {\it length} of a simple diagram is $n = \sum n_{ij}, 1 \le i
< j \le m$.

The Birman--Series parameterization of elements of $J_0$ consists
of two sets of data. The first is the ordered sequence
$h_1(\gamma)=(n_{12}, n_{13}, \cdots, n_{m-1,m})$ which records
for each pair of distinct elements $a_i, a_j$ of $R$ the number
$n_{ij}$ of segments of $\gamma$ which join $a_i$ to $a_j$. The
second set of data, $h_2(\gamma)$, records information about the
position of the initial and final points $t_0, t_n$ of $\gamma$.
Let $a(t_i)$ be the element of $A$ containing $t_i$ and let
$j(t_i) \in \mathbf N$ be the position of $t_i$ among the
partition points of $\gamma$ which lie along $a(t_i)$ counting in
the anticlockwise direction round $\partial R$. Define
$h_2(\gamma)=(a(t_0), j(t_0), a(t_n), j(t_n))$.

The following lemmas and their proofs in \cite{birman-series1985t}
still hold in our case.

\begin{lem} Suppose that $\gamma, \gamma^{\prime} \in J_0$ and
that $h_1(\gamma)=h_1(\gamma^{\prime}),
h_2(\gamma)=h_2(\gamma^{\prime})$. Let $t_0, t_1, \cdots, t_n$ and
$t_0^{\prime}, t_1^{\prime}, \cdots, t_n^{\prime}$ be the
partition points of $\gamma, \gamma^{\prime}$ respectively. Then
$a(t_i)=a(t_i^{\prime})$ for each $i=0,1,\cdots,n$.
\end{lem}

\begin{lem} Let $J_0(n)=\{\gamma \in J_0 : \, \parallel \gamma \parallel =n
\}$. Then there is a polynomial $P_0(n)$ such that the number of
simple diagrams of length $n$
$$
{\rm card} \{(h_1(\gamma), h_2(\gamma)): \gamma \in J_0(n)\} \le
P_0(n).
$$
\end{lem}

The main idea of the proof of Birman--Series Theorem in
\cite{birman-series1985t} is that geodeisc arcs in $J_0(n)$ (for
sufficiently large $n$) with the same parameterization lie
exponentially close in $M$. It relies on the following key lemma
which is Lemma 3.1 in \cite{birman-series1985t}.

\begin{lem} There is a universal constant $\alpha >0$ (depending
only on the choice of the fundamental polygon $R$) so that
$$
l(\gamma) \ge \alpha \parallel \gamma \parallel
$$
for $\gamma \in J_0$ with $\parallel \gamma \parallel$
sufficiently large, where $l(\gamma)$ denotes the hyperbolic
length of $\gamma$.
\end{lem}

\begin{pf}
There is a universal constant $\epsilon >0$ so that any segment of
$\gamma$ which does not connect two consecutive sides of $R$ or
does not intersect a suitably chosen disk neighborhood of each
cusp or cone point has hyperbolic length at least $\epsilon$. Let
$q$ be the maximum number of sides of $R$, projected to $M$, which
meet at any cusp or cone point of $M$. Then at most $q-1$
consecutive segments of $\gamma$ can connect consecutive sides of
$R$ around the same cusp or cone point and intersect the chosen
disk neighborhood of that cusp or cone point; for otherwise there
will be a self-intersection on $\gamma$. Hence in any $q$
consecutive segments of $\gamma$, at least one has hyperbolic
length $\epsilon$, which gives the result.
\end{pf}

The following two lemmas then apply respectively to the set of all
complete simple geodesics which never intersect any cusp or cone
point and to the set of simple geodesics which start from a fixed
cusp or cone point and never terminates at any cusp or cone point.
(Recall that we assume that $M$ has no boundary geodesics.)

\begin{lem}\label{lem:BS 3.2}
Let $\gamma, \gamma^{\prime} \in J_0(2n+1)$ and suppose that
$h_1(\gamma)=h_1(\gamma^{\prime}),
h_2(\gamma)=h_2(\gamma^{\prime})$. Let $\delta \subset \gamma,
\delta^{\prime} \subset \gamma^{\prime}$ denote the segments of
$\gamma, \gamma^{\prime}$ lying between the partition points $t_n,
t_{n+1}$ and $t_n^{\prime}, t_{n+1}^{\prime}$ respectively. Then
$\delta^{\prime} \subset B_{ce^{-\alpha n}}(\delta)$ where $c,
\alpha$ are universal constants and where $B_\epsilon(\delta)$
denotes the tubular neighborhood of $\delta$ of hyperbolic radius
$\epsilon>0$.
\end{lem}

\begin{lem}\label{lem:gBS 3.2}
Let $\gamma, \gamma^{\prime} \in J_0(n+k)$ be such that they start
at the same vertex of $R$ and that
$h_1(\gamma)=h_1(\gamma^{\prime}),
h_2(\gamma)=h_2(\gamma^{\prime})$. Let $\delta \subset \gamma,
\delta^{\prime} \subset \gamma^{\prime}$ denote the segments of
$\gamma, \gamma^{\prime}$ lying between the partition points $t_i,
t_{i+1}$ and $t_i^{\prime}, t_{i+1}^{\prime}$ respectively for
some $1 \le i \le k$. Then $\delta^{\prime} \subset B_{ce^{-\alpha
n}}(\delta)$ where $c, \alpha$ are universal constants and where
$B_\epsilon(\delta)$ denotes the tubular neighborhood of $\delta$
of hyperbolic radius $\epsilon>0$.
\end{lem}

\noindent Note that Lemma \ref{lem:BS 3.2} is Lemma 3.2 in
\cite{birman-series1985t} and Lemma \ref{lem:gBS 3.2} can be
proved similarly.

\vskip 10pt

From these we have the following proposition which is Proposition
4.1 in \cite{birman-series1985t} from which the conclusion of the
Birman--Series Theorem follows exactly as in the proofs in
\cite{birman-series1985t} \S 5.

\begin{prop}
There exist universal constants $L, c, \alpha >0$ and a polynomial
$P_0(\cdot)$ such that for each $n$ there is a set $F_n$ of simple
geodesic arcs, each of length at most $L$, so that ${\rm
card}(F_n) \le P_0(n)$ and so that
$$
S_0 \subset \cup \{B_\epsilon(\gamma) \mid \gamma \in F_n\},
\epsilon = ce^{\alpha n}.
$$
\end{prop}

\vskip 15pt

Finally we remark that the above Birman--Series' arguments will
give rough estimates on the distribution of simple closed
geodesics on a compact hyperbolic cone-surface $M$ which is enough
for proving the absolute convergence of the series appearing in
various generalized McShane's identities, as was observed and used
in \cite{akiyoshi-miyachi-sakuma2004preprint} (for the case of
complete hyperbolic surfaces) for similar purposes.

\vskip 10pt

\begin{lem}\label{lem:BS}
Let $M$ be a compact hyperbolic cone-surface with all cone angles
in $(0,\pi]$. Then for any constant $c>0$
\vskip 6pt

{\rm (i)} the series
\begin{eqnarray*}
\sum_{\beta} \frac{1}{\exp(c|\beta|)}
\end{eqnarray*}
converges absolutely, where the sum is over all generalized simple
closed geodesics on $M$ and all simple normal geodesic arcs
connecting geometric boundary components of $M$;

\vskip 6pt
{\rm (ii)} the series
\begin{eqnarray*}
\sum_{\alpha, \beta} \frac{1}{\exp[c(|\alpha|+|\beta|)]}
\end{eqnarray*}
converges absolutely, where the sum is over all pairs $\alpha,
\beta$ of disjoint generalized simple closed geodesics on $M$
and/or simple normal geodesic arcs connecting geometric boundary
components of $M$.\square
\end{lem}

\vskip 30pt
\section{{\bf Proof of Theorems}}\label{s:proof}
\vskip 30pt

{\it Proof of Theorem \ref{thm:mcshane most general}} \,\, Now the
proof is obvious from the previous discussions. Suppose $\Delta_0$
is a cone point. Recall $\mathcal H$ is a suitably chosen small
circle centered at $\Delta_0$, and ${\mathcal H}_{\texttt{ns}}$,
${\mathcal H}_{\texttt{sn}}$, ${\mathcal H}_{\texttt{snn}}$ are
the point sets of the first intersections of $\mathcal H$ with
respectively all non-simple, all simple-normal, all
simple-not-normal $\Delta_{0}$-geodesics. The elliptic measure of
each of these subsets of $\mathcal H$ is the radian measure that
it subtends to the cone point $\Delta_0$. The generalized
Birman--Series Theorem in \S \ref{s:gBS} implies that the closed
subset ${\mathcal H}_{\texttt{sn}}$ has measure $0$. Hence the
open subset ${\mathcal H}_{\texttt{ns}} \cup {\mathcal
H}_{\texttt{snn}}$ has full measure, that is, $\theta_0$. Now the
maximal open intervals of ${\mathcal H}_{\texttt{ns}} \cup
{\mathcal H}_{\texttt{snn}}$, suitably combined, have measure $2
{\rm Gap}(\Delta_{0};\alpha, \beta)$ for each unordered pair of
generalized simple closed geodesics $\alpha, \beta$ on $M$ which
bound with $\Delta_0$ an embedded pair of pants on $M$. Hence
their sum is equal to $\theta_0$ and the desired identity follows.
The cases where $\Delta_0$ is a boundary geodesic or a cusp are
similarly proved.  \square

\vskip 12pt

{\it Proof of Corollary \ref{cor:mcshane conical holed
weierstrass}} \,\, Consider the case where $\Delta_0$ is a cone
point. In this case $T$ admits a unique elliptic involution $\eta$
such that $\eta$ maps each oriented simple closed geodesics on $T$
onto itself with orientation reversed. Note that $\eta$ fixes the
cone point $\Delta_0$ and three other interior points which are
the so-called Weierstrass points of $T$. Each simple closed
geodesics on $T$ passes exactly two Weierstrass points; hence
there are three Weierstrass classes of simple closed geodesics on
$T$. Now the quotient of $T$ under $\eta$ is a sphere with three
angle $\pi$ cone points and a cone point with angle $\theta /2$.
Then Theorem \ref{thm:mcshane most general} applies to $M =
T/\langle\eta\rangle$, with $\Delta_0$ the angle $\pi$ cone point
whose inverse image under $\eta$ is the Weierstrass point that the
Weierstrass class $\mathcal A$ misses. Note that each generalized
simple closed geodesic on $M = T/\langle\eta\rangle$ is either a
geometric boundary component or degenerate simple closed geodesic
which is the double cover of a simple geodesic arc which connects
two Weierstrass points. Hence the set of all pairs of generalized
simple closed geodesics which bound with $\Delta_0$ an embedded
pair of pants is exactly the set of pairs consisting of the angle
$\theta /2$ cone point plus a degenerate simple closed geodesic
$\gamma^{\prime}$ which is the double cover of the quotient simple
geodesic arc of a simple closed geodesic $\gamma$ on $T$ in the
given Weierstrass class $\mathcal A$ (note that by definition the
length of $\gamma^{\prime}$ is the same as that of $\gamma$).
Hence by (\ref{subcase 1.3}) the summand in the summation is
\begin{eqnarray*}
\frac{\pi}{2} - \tan^{-1} \left (
\frac{\sin\frac{\pi}{2}\sinh\frac{|\gamma|}{2}
}{\cos\frac{\theta}{4}+\cos\frac{\pi}{2}\cosh\frac{|\gamma|}{2} }
\right ) = \tan^{-1} \left (
\frac{\cos\frac{\theta}{4}}{\sinh\frac{|\gamma|}{2}} \right ).
\end{eqnarray*}

The proof for the case where $\Delta_0$ is a boundary geodesic is
similar. \square

\begin{rmk}
Note that we can also choose $\Delta_0$ to be the angle $\theta
/2$ cone point on $T/\langle\eta\rangle$, then we obtain
(\ref{eqn:mcshane cone torus}), the generalization of McShane's
original identity to the cone-torus $T$. This is one way of seeing
why we can allow the cone angle of up to $2\pi$ in the cone torus
case.
\end{rmk}

\vskip 10pt

{\it Proof of Theorem \ref{thm:mcshane genus two global}} \,\, It
is well known that $M$ admits a unique hyperelliptic involution
$\eta$ (see for example \cite{haas-susskind1989pams}) such that
$\eta$ maps each simple closed geodesic onto itself and
preserves/reverses the orientation of separating/non-separating
simple closed geodesics. Note that $\eta$ leaves six points on $M$
fixed; they are the six Weierstrass points on $M$. Consider the
quotient $M^{\prime}=M/\langle\eta\rangle$ which is a sphere with
six angle $\pi$ cone points. Each generalized simple closed
geodesic on $M^{\prime}$ is either \begin{itemize} \item[(i)] an
angle $\pi$ cone point; or \item[(ii)] a degenerate simple closed
geodesic ${\beta}^{\prime}$ which is the double cover of a simple
geodesic arc $c$ connecting two angle $\pi$ cone points where the
inverse image of $c$ under $\eta$ is a non-separating simple
closed geodesic $\beta$ on $M$; or
\item[(iii)] a separating (non-degenerate) simple closed geodesic
${\alpha}^{\prime}$ whose inverse image under $\eta$ is a
separating simple closed geodesic $\alpha$ on $M$. In this case
${\alpha}^{\prime}$ does not pass through any of the six angle
$\pi$ cone points and there are three of them on each side of
${\alpha}^{\prime}$ on $M^{\prime}$. Hence $\alpha$ passes none of
six Weierstrass points and there are three of them on each side of
$\alpha$ on $M$.
\end{itemize}
Now apply Theorem \ref{thm:mcshane most general} to $M^{\prime}$
with $\Delta_0$ one of the six angle $\pi$ cone points. Then each
pair of generalized simple closed geodesics on $M^{\prime}$ which
bound with $\Delta_0$ an embedded pair of pants $\mathcal P$
consists of a separating simple closed geodesic $\alpha^{\prime}$
on $M^{\prime}$ and a degenerate simple closed geodesic
$\beta^{\prime}$ on $M^{\prime}$ which lies on the same side of
$\alpha^{\prime}$ as $\Delta_0$ and misses $\Delta_0$. Let the
inverse image of $\alpha^{\prime}, \beta^{\prime}$ under $\eta$ be
$\alpha, \beta$. Then $\alpha$ is a separating simple closed
geodesic on $M$ and $\beta$ is a non-separating simple closed
geodesic on $M$. Furthermore, $\beta$ and the Weierstrass point
which is the inverse image of $\Delta_0$ lie on the same side of
$\alpha$ on $M$. Note that the hyperbolic lengths of
$\alpha^{\prime}, \beta^{\prime}$ are respectively $|\alpha|/2,
|\beta|$. Hence by (\ref{eqn:dGf 1.1}) in this case the summand in
the resulting generalized McShane's Weierstrass identity for
$M^{\prime}$ with the chosen $\Delta_0$ is
\begin{eqnarray*}
2 \tan^{-1} \left (\frac{\sin\frac{\pi}{2}}{\cos\frac{\pi}{2}+
\exp\frac{|\alpha|/2+|\beta|}{2}} \right ) = 2 \tan^{-1}
\exp\left(-\frac{|\alpha|}{4}-\frac{|\beta|}{2}\right).
\end{eqnarray*}
Note that each pair of disjoint simple closed geodesics
$(\alpha,\beta)$ on $M$ such that $\alpha$ is separating and
$\beta$ is non-separating arises as the inverse image of a unique
pair of generalized simple closed geodesics on $M'$ as described
above, where the chosen $\Delta_0$ is the angle $\pi$ cone point
which is the image under $\eta$ of the Weierstrass point on $M$
that lies on the same side of $\alpha$ as $\beta$ and is missed by
$\beta$.

Summing all the six resulting Weierstrass identities we then have
\begin{eqnarray*}
\sum 2 \tan^{-1} \exp \left( -\frac{|\alpha|}{4}-
\frac{|\beta|}{2} \right ) = \frac{6\pi}{2},
\end{eqnarray*}
where the sum is over all ordered pairs $(\alpha, \beta)$ of
disjoint simple closed geodesics on $M$ such that $\alpha$ is
separating and $\beta$ is non-separating. \square

\vskip 15pt

{\it Proof of Addendum \ref{add:genus two}} \,\, We first prove
that the series in (\ref{eqn:genus two QF}) converges absolutely
and uniformly on compact set in the space $\mathcal Q \mathcal F$
of quasi-Fuchsian representations of ${\pi_1}(M)$ into ${\rm
SL}(2, \mathbf C)$ by the same argument as used in
\cite{akiyoshi-miyachi-sakuma2004preprint}. The identity
(\ref{eqn:genus two QF}) then follows by analytic continuation
since each summand in it is an analytic function of the complex
Fenchel--Nielsen coordinates for the quasi-Fuchsian space (see
\cite{tan1994ijm}) and the identity holds when all the coordinates
take real values (by Theorem \ref{thm:mcshane genus two global})
and the space of quasi-Fuchsian representations of ${\pi_1}(M)$
into ${\rm PSL}(2, \mathbf C)$ is simply connected.

As pointed out in \cite{akiyoshi-miyachi-sakuma2004preprint} Lemma
5.2, by \cite{jorgensen-marden1979qjm} Lemma 3, for any compact
subset $\mathcal C$ of $\mathcal Q \mathcal F$, there is a
constant $k = k(C) >0$ such that
\begin{eqnarray*}
k l_{\rho_0}(\gamma) \le \Re l_{\rho}(\gamma) \le k^{-1}
l_{\rho_0}(\gamma),
\end{eqnarray*}
for any essential simple closed curve $\gamma$, where $\rho_0$ is
a fixed Fuchsian representation of ${\pi_1}(M)$ into ${\rm SL}(2,
\mathbf C)$.

Since $|\tan^{-1}(x)| \le 2|x|$ for $|x|$ sufficiently small, we
have for all except a finitely many pairs of (free homotopy
classes of) disjoint essential simple closed curves $\alpha,
\beta$ on $M$ such that $\alpha$ is separating and $\beta$ is
non-separating
\begin{eqnarray*}
\left|
\tan^{-1}\exp\Big(-\frac{l_{\rho}(\alpha)}{4}-\frac{l_{\rho}(\beta)}{2}\Big)\right|
&\le&
2\left|\exp\Big(-\frac{l_{\rho}(\alpha)}{4}-\frac{l_{\rho}(\beta)}{2}\Big)\right|\\
&=& 2 \exp\Big(-\frac{\Re l_{\rho}(\alpha)}{4}-\frac{\Re
l_{\rho}(\beta)}{2}\Big)\\ &\le& 2
\exp\Big(-k\Big(\frac{l_{\rho_0}(\alpha)}{4}+\frac{l_{\rho_0}(\beta)}{2}\Big)\Big).
\end{eqnarray*}
Thus the series in (\ref{eqn:genus two QF}) converges absolutely
and uniformly on the compact set $C$ of $\mathcal Q \mathcal F$
since the series
$$
\sum
\exp\Big(-k\Big(\frac{l_{\rho_0}(\alpha)}{4}+\frac{l_{\rho_0}(\beta)}{2}\Big)\Big)
$$
converges by Lemma \ref{lem:BS}. \square

\vskip 10pt

\vskip 30pt
\section{{\bf Complexified reformulation of the generalized McShane's
identity }}\label{s:reformulation} \vskip 30pt

In this section we prove the unified version (\ref{eqn:reform of
cp and gb cases}) of our generalized McShane's identity using
complex arguments and interpret it geometrically.

\vskip 8pt

\noindent {\bf Two functions} \,\, First we would like to define
two functions $G, S: {\mathbf C}^3 \rightarrow \mathbf C$ as
follows:
\begin{eqnarray}
G(x,y,z)=2
\tanh^{-1}\left(\frac{\sinh(x)}{\cosh(x)+\exp(y+z)}\right),
\end{eqnarray}
\begin{eqnarray}
S(x,y,z)=\tanh^{-1}\left(\frac{\sinh(x)\sinh(y)}{\cosh(z)+\cosh(x)\cosh(y)}\right).
\end{eqnarray}
Note that here for a complex number $x$, $\tanh^{-1}(x)$ is
defined to have imaginary part in $(-\pi/2, \pi/2]$. Using the
identity
$$
x=\frac{1}{2}\log\frac{1+\tanh(x)}{1-\tanh(x)},
$$
it is easy to check that the two functions have also the following
expressions:
\begin{eqnarray}
G(x,y,z)=\log\frac{\exp(x)+\exp(y+z)}{\exp(-x)+\exp(y+z)},
\end{eqnarray}
\begin{eqnarray}
S(x,y,z)=\frac{1}{2}\log\frac{\cosh(z)+\cosh(x+y)}{\cosh(z)+\cosh(x-y)},
\end{eqnarray}
as used by Mirzakhani in \cite{mirzakhani2004preprint}. (She uses
different notations $\mathcal D, \mathcal R$ as explained below.)
Here for a non-zero complex number $x$, $\log(x)$ assumes the main
branch value with imaginary part in $(-\pi, \pi]$. We shall see
that both expressions of the functions are useful.

\vskip 10pt

For $x,y,z >0$, the geometrical meanings of $G(x,y,z)$ and
$S(x,y,z)$ are as follows. Let $\mathcal P(2x,2y,2z)$ be the
unique hyperbolic pair of pants whose boundary components $X,Y,Z$
are simple closed geodesics of lengths $2x,2y,2z$ respectively.
Then $S(x,y,z)$ is half the length of the orthogonal projection of
the boundary geodesic $Y$ onto $X$ in $\mathcal P(2x,2y,2z)$ and
$S(x,z,y)$ is half the length of the orthogonal projection of the
boundary geodesic $Z$ onto $X$ in $\mathcal P(2x,2y,2z)$, and
$G(x,y,z)$ is the length of each of the two gaps between these two
projections on $X$. We have therefore the identity
\begin{eqnarray}
G(x,y,z)+S(x,y,z)+S(x,z,y)=x
\end{eqnarray}
for all $x,y,z \ge 0$. Note that the same identity holds modulo
$\pi i$ for all $x,y,z \in \mathbf C$.

\vskip 10pt

\begin{rmk} The relations of our functions $G,S$ with Mirzakhani's
functions $\mathcal D, \mathcal R$ are
\begin{eqnarray}
G(x,y,z)=\mathcal D(2x,2y,2z) /2,
\end{eqnarray}
\begin{eqnarray}
S(x,y,z)=(x-\mathcal R(2x,2z,2y)) /2.
\end{eqnarray}
\end{rmk}

\vskip 10pt

\begin{lem}\label{lem: Gap=G+S}
{\rm (i)} For $x,z \ge 0$ and $y \in [0, \frac{\pi}{2}]$,
\begin{eqnarray}\label{eqn:(x,yi,z)}
\hskip 6pt G(x,yi,z)+S(x,yi,z)=x-\tanh^{-1}\left(
\frac{\sinh(x)\sinh(z)}{\cos(y)+\cosh(x)\cosh(z)} \right).
\end{eqnarray}
{\rm (ii)} For $x,y \in [0, \frac{\pi}{2}]$ and $z \ge 0$,
\begin{eqnarray}\label{eqn:(xi,yi,z)}
\hskip 16pt G(xi,yi,z)+S(xi,yi,z)=\left[x-\tan^{-1}\left(
\frac{\sin(x)\sinh(z)}{\cos(y)+\cos(x)\cosh(z)} \right)\right]i.
\end{eqnarray}
\end{lem}

\vskip 10pt

\begin{pf} {\rm (i)} It follows from the following two identities since $\Re S(x,yi,z)=0$:
\begin{eqnarray}\label{eqn:(i)1}
\Re G(x,yi,z)=x-\tanh^{-1}\left(
\frac{\sinh(x)\sinh(z)}{\cos(y)+\cosh(x)\cosh(z)} \right),
\end{eqnarray}

\begin{eqnarray}\label{eqn:(i)2}
\Im G(x,yi,z)+ \Im S(x,yi,z)=0.
\end{eqnarray}

\vskip 8pt

{\it Proof of {\rm (\ref{eqn:(i)1})} and {\rm (\ref{eqn:(i)2})}}:
\,\, By definition,
\begin{eqnarray*}
G(x,yi,z)&=&\log\frac{\exp(x)+\exp(yi+z)}{\exp(-x)+\exp(yi+z)}\\
&=&\log\frac{[\exp(x)+\cos(y)\exp(z)]+i[\sin(y)\exp(z)]}{[\exp(-x)+\cos(y)\exp(z)]+i[\sin(y)\exp(z)]}.
\end{eqnarray*}
Hence
\begin{eqnarray*}
\Re G(x,yi,z)&=&\frac{1}{2}\log\frac{[\exp(x)+\cos(y)\exp(z)]^2+
[\sin(y)\exp(z)]^2}{[\exp(-x)+\cos(y)\exp(z)]^2+[\sin(y)\exp(z)]^2}\\
&=&\frac{1}{2}\log\frac{\exp(2x)+\exp(2z)+2\exp(x)\cos(y)\exp(z)}{\exp(-2x)+\exp(2z)+2\exp(-x)\cos(y)\exp(z)}\\
&=&\frac{1}{2}\log\left(\frac{\cosh(x-z)+\cos(y)}{\cosh(x+z)+\cos(y)}\frac{\exp(x+z)}{\exp(-x+z)}\right)\\
&=&x-\frac{1}{2}\log\frac{\cosh(x+z)+\cos(y)}{\cosh(x-z)+\cos(y)}\\
&=&x-\tanh^{-1}\left(\frac{\sinh(x)\sinh(z)}{\cos(y)+\cosh(x)\cosh(z)}\right).
\end{eqnarray*}

On the other hand,
\begin{eqnarray*}
\Im G(x,yi,z) &=&
\tan^{-1}\left(\frac{\sin(y)\exp(z)}{\exp(x)+\cos(y)\exp(z)}\right)
- \tan^{-1}\left(\frac{\sin(y)\exp(z)}{\exp(-x)+\cos(y)\exp(z)}\right)\\
&=& \tan^{-1}\left(\frac{[\exp(-x)-\exp(x)]\sin(y)\exp(z)}{[\exp(x)+\cos(y)\exp(z)][\exp(-x)+\cos(y)\exp(z)]+[\sin(y)\exp(z)]^2}\right)\\
&=& \tan^{-1}\left(\frac{[\exp(-x)-\exp(x)]\sin(y)\exp(z)}{1+\exp(2z)+[\exp(x)+\exp(-x)]\cos(y)\exp(z)}\right)\\
&=& - \tan^{-1}\left(\frac{\sinh(x)\sin(y)}{\cosh(z)+\cosh(x)\cos(y)}\right)\\
&=& -\Im S(x,yi,z),
\end{eqnarray*}
since
\begin{eqnarray*}
S(x,yi,z)&=&\tanh^{-1}\left(\frac{\sinh(x)\sinh(yi)}{\cosh(z)+\cosh(x)\cosh(yi)}\right)\\
&=& i
\,\tan^{-1}\left(\frac{\sinh(x)\sin(y)}{\cosh(z)+\cosh(x)\cos(y)}\right).
\end{eqnarray*}

\vskip 10pt

{\rm (ii)} It will follow from the following two identities:
\begin{eqnarray}\label{eqn:(ii)1}
\Im G(xi,yi,z)=x-\tan^{-1}\left(
\frac{\sin(x)\sinh(z)}{\cos(y)+\cos(x)\cosh(z)} \right),
\end{eqnarray}

\begin{eqnarray}\label{eqn:(ii)2}
\Re G(xi,yi,z)+ S(xi,yi,z)=0.
\end{eqnarray}

\vskip 8pt

{\it Proof of {\rm (\ref{eqn:(ii)1})} and {\rm
(\ref{eqn:(ii)2})}}:  \,\, By definition,
\begin{eqnarray*}
G(xi,yi,z)&=&\log\frac{\exp(xi)+\exp(yi+z)}{\exp(-xi)+\exp(yi+z)}\\
&=&\log\frac{[\cos(x)+\cos(y)\exp(z)]+i[\sin(x)+\sin(y)\exp(z)]}{[\cos(x)+\cos(y)\exp(z)]+i[-\sin(x)+\sin(y)\exp(z)]}.
\end{eqnarray*}
Hence
\begin{eqnarray*}
\Re G(xi,yi,z)&=&\frac{1}{2}\log\frac{[\cos(x)+\cos(y)\exp(z)]^2+
[\sin(x)+\sin(y)\exp(z)]^2}{[\cos(x)+\cos(y)\exp(z)]^2+[-\sin(x)+\sin(y)\exp(z)]^2}\\
&=&\frac{1}{2}\log\frac{1+\exp(2z)+\cos(x-y)\exp(z)}{1+\exp(2z)+\cos(x+y)\exp(z)}\\
&=&\frac{1}{2}\log\frac{\cosh(z)+\cos(x-y)}{\cosh(z)+\cos(x+y)}\\
&=&-\frac{1}{2}\log\frac{\cosh(z)+\cosh(xi+yi)}{\cosh(z)+\cosh(xi-yi)}\\
&=&-S(xi,yi,z).
\end{eqnarray*}

On the other hand,
\begin{eqnarray*}
I&=&\Im G(xi,yi,z)\\
&=&\tan^{-1}\left(\frac{\sin(x)+\sin(y)\exp(z)}{\cos(x)+\cos(y)\exp(z)}\right)
- \tan^{-1}\left(\frac{-\sin(x)+\sin(y)\exp(z)}{\cos(x)+\cos(y)\exp(z)}\right)\\
&=& \tan^{-1}\left(\frac{2\sin(x)[\cos(x)+\cos(y)\exp(z)]}{[\cos(x)+\cos(y)\exp(z)]^2-[\sin(x)]^2+[\sin(y)\exp(z)]^2}\right)\\
&=&\tan^{-1}\left(\frac{\sin(2x)+2\sin(x)\cos(y)\exp(z)}{\cos(2x)+\exp(2z)+2\cos(x)\cos(y)\exp(z)}\right).
\end{eqnarray*}
Hence
\begin{eqnarray*}
iI=\tanh^{-1}\left(\frac{i\sin(2x)+2i\sin(x)\cos(y)\exp(z)}{\cos(2x)+\exp(2z)+2\cos(x)\cos(y)\exp(z)}\right),
\end{eqnarray*} or
\begin{eqnarray*}
\frac{\exp(2iI)-1}{\exp(2iI)+1}=\frac{i\sin(2x)+2i\sin(x)\cos(y)\exp(z)}{\cos(2x)+\exp(2z)+2\cos(x)\cos(y)\exp(z)}.
\end{eqnarray*}
Hence
\begin{eqnarray*}
\exp(2iI)&=&\frac{\exp(2xi)+\exp(2z)+2\exp(xi)\cos(y)\exp(z)}{\exp(-2xi)+\exp(2z)+2\exp(-xi)\cos(y)\exp(z)}\\
&=&\frac{\cosh(xi-z)+\cos(y)}{\cosh(xi+z)+\cos(y)}\frac{\exp(xi+z)}{\exp(-xi+z)}\\
&=&\exp(2xi)\frac{\cosh(yi)+\cosh(xi-z)}{\cosh(yi)+\cosh(xi+z)}.
\end{eqnarray*}
Thus
\begin{eqnarray*}
iI&=&xi-\frac{1}{2}\log\frac{\cosh(yi)+\cosh(xi+z)}{\cosh(yi)+\cosh(xi-z)}\\
&=&xi-\tanh^{-1}\left(\frac{\sinh(xi)\sinh(z)}{\cosh(yi)+\cosh(xi)\cosh(z)} \right)\\
&=&xi-i\tan^{-1}\left(\frac{\sin(x)\sinh(z)}{\cos(y)+\cos(x)\cosh(z)}\right).
\end{eqnarray*}
Therefore
\begin{eqnarray*}
\Im G(xi,yi,z)=I=
x-\tan^{-1}\left(\frac{\sin(x)\sinh(z)}{\cos(y)+\cos(x)\cosh(z)}\right).
\end{eqnarray*}
\end{pf}

\vskip 10pt

\noindent {\bf Restatement of the complexified identities} \,\,
Now we can restate the non-cusp cases of Theorem
\ref{thm:complexified} using the functions $G,S$ defined above.
Recall that for each generalized simple closed geodesic $\delta$,
we have defined in \S \ref{s:intro} its {\bf complex length}
$|\delta|$, that is, $|\delta|=0$ if $\delta$ is a cusp;
$|\delta|=\theta i$ if $\delta$ is a cone point of angle $\theta
\in (0, \pi]$; and $|\delta|=l$ if $\delta$ is a boundary geodesic
or an interior generalized simple closed geodesic of hyperbolic
length $l>0$.

\vskip 10pt

\begin{thm}\label{thm:complexified non-cusp cases}
For a compact hyperbolic cone-surface $M$ with all cone angles in
$(0, \pi]$, let all its geometric boundary components be
$\Delta_0, \Delta_1, \cdots, \Delta_n$ with complex lengths $L_0,
L_1, \cdots, L_n$ respectively. If $\Delta_0$ is a cone point or a
boundary geodesic then
\begin{eqnarray}\label{eqn:reform of non-cusp cases with GS}
\sum_{\alpha, \beta}
G\left(\frac{L_0}{2},\frac{|\alpha|}{2},\frac{|\beta|}{2}\right) +
\sum_{j=1}^{n}\sum_{\beta}
S\left(\frac{L_0}{2},\frac{L_j}{2},\frac{|\beta|}{2}\right)=
\frac{L_0}{2},
\end{eqnarray} where the first sum is over all
(unordered) pairs of generalized simple closed geodesics $\alpha,
\beta$ on $M$ such that $\alpha, \beta$ bound with $\Delta_0$ an
embedded pair of pants on $M$ (note that one of $\alpha, \beta$
might be a geometric boundary component) and the sub-sum in the
second sum is over all interior simple closed geodesics $\beta$
such that $\beta$ bounds with $\Delta_j$ and $\Delta_0$ an
embedded pair of pants on $M$. Furthermore, all the series in
(\ref{eqn:reform of non-cusp cases with GS}) converge absolutely.
\end{thm}

\vskip 10pt

\begin{rmk}
We shall omit the proof of Theorem \ref{thm:complexified} in the
case where $\Delta_0$ is a cusp, for as remarked before, in the
cusp case the identity (\ref{eqn:reform of cusp cases}) can either
be proved similarly or be derived by considering the first order
infinitesimal terms of the corresponding identity (\ref{eqn:reform
of cp and gb cases}) in other cases.
\end{rmk}

\vskip 10pt

\begin{pf}
We first show that our generalized McShane's identities
(\ref{eqn:001}) and (\ref{eqn:002}) can be reformulated as
(\ref{eqn:reform of non-cusp cases with GS}) modulo convergence.

\vskip 8pt

First suppose that $\Delta_0$ is a boundary geodesic of hyperbolic
length $l_0 >0$.

For a pair of interior generalized simple closed geodesics
$\alpha, \beta$ which bound with $\Delta_0$ an embedded pair of
pants on $M$, we have directly by definition that
\begin{eqnarray*}
{\rm Gap}(\Delta_0;\alpha,\beta)=
G\left(\frac{l_0}{2},\frac{|\alpha|}{2},\frac{|\beta|}{2}\right).
\end{eqnarray*}

For a pair of generalized simple closed geodesics $\alpha, \beta$
such that $\alpha$ is a boundary geodesic and $\beta$ is an
interior generalized simple closed geodesic and that they bound
with $\Delta_0$ an embedded pair of pants on $M$, we have by
definition and the geometric meanings of $G,S$ that
\begin{eqnarray*}
{\rm Gap}(\Delta_0;\alpha,\beta)=
G\left(\frac{l_0}{2},\frac{|\alpha|}{2},\frac{|\beta|}{2}\right)+
S\left(\frac{l_0}{2},\frac{|\alpha|}{2},\frac{|\beta|}{2}\right).
\end{eqnarray*}

For a pair of generalized simple closed geodesics $\alpha, \beta$
such that $\alpha$ is a cone point of angle $\varphi \in (0,\pi]$
and $\beta$ is an interior generalized simple closed geodesic and
that they bound with $\Delta_0$ an embedded pair of pants on $M$,
we have by (\ref{eqn:(x,yi,z)}) with $x=l_0/2, y=\varphi/2,
z=|\beta|/2$ that
\begin{eqnarray*}
{\rm Gap}(\Delta_0;\alpha,\beta)=
G\left(\frac{l_0}{2},\frac{\varphi i}{2},\frac{|\beta|}{2}\right)+
S\left(\frac{l_0}{2},\frac{\varphi i}{2},\frac{|\beta|}{2}\right).
\end{eqnarray*}

\vskip 10pt

Next suppose that $\Delta_0$ is a cone point of angle $\theta_0
\in (0,\pi]$.

\vskip 8pt

For a pair of interior generalized simple closed geodesics
$\alpha, \beta$ which bound with $\Delta_0$ an embedded pair of
pants on $M$, we have by definition that
\begin{eqnarray*}
{\rm Gap}(\Delta_0;\alpha,\beta)\,i=
G\left(\frac{\theta_0i}{2},\frac{|\alpha|}{2},\frac{|\beta|}{2}\right).
\end{eqnarray*}

For a pair of generalized simple closed geodesics $\alpha, \beta$
such that $\alpha$ is a boundary geodesic and $\beta$ is an
interior generalized simple closed geodesic and that they bound
with $\Delta_0$ an embedded pair of pants on $M$, we have by the
analysis in \S \ref{s:calculation} that
\begin{eqnarray*}
& &{\rm Gap}(\Delta_0;\alpha,\beta)\,i\\&=&2i\tan^{-1}
\frac{\sin\frac{\theta_0}{2}}{\cos\frac{\theta_0}{2}+\exp\frac{|\alpha|+|\beta|}{2}}
+i\tan^{-1} \frac{\sin\frac{\theta_0}{2}\sinh\frac{|\alpha|}{2}}
{\cosh\frac{|\beta|}{2}+\cos\frac{\theta_0}{2}\cosh\frac{|\alpha|}{2}}\\
&=&G\left(\frac{\theta_0i}{2},\frac{|\alpha|}{2},\frac{|\beta|}{2}\right)+
   S\left(\frac{\theta_0i}{2},\frac{|\alpha|}{2},\frac{|\beta|}{2}\right).
\end{eqnarray*}

For a pair of generalized simple closed geodesics $\alpha, \beta$
such that $\alpha$ is a cone point of angle $\varphi \in (0,\pi]$
and $\beta$ is an interior generalized simple closed geodesic and
that they bound with $\Delta_0$ an embedded pair of pants on $M$,
we have by (\ref{eqn:(xi,yi,z)}) with $x=\theta_0/2, y=\varphi/2,
z=|\beta|/2$ that
\begin{eqnarray*}
{\rm Gap}(\Delta_0;\alpha,\beta)\,i =
G\left(\frac{\theta_0i}{2},\frac{\varphi i}{2},\frac{|\beta|}{2}
\right)+ S\left(\frac{\theta_0i}{2},\frac{\varphi i}{2},
\frac{|\beta|}{2} \right).
\end{eqnarray*}

\vskip 10pt

Finally we prove the absolute convergence of the series in
(\ref{eqn:reform of non-cusp cases with GS}). It is not hard to
see that we only need to prove, for each $j=1, \cdots, n$, the
absolute convergence of the series
\begin{eqnarray*}
\sum_{\beta}
S\left(\frac{L_0}{2},\frac{L_j}{2},\frac{|\beta|}{2}\right),
\end{eqnarray*}
where the sum is over all interior generalized simple closed
geodesics $\beta$ which bounds with $\Delta_j$ and $\Delta_0$ an
embedded pair of pants on $M$. The desired absolute convergence
follows from Lemma \ref{lem:BS} since
\begin{eqnarray*}
S\left(\frac{L_0}{2},\frac{L_j}{2},\frac{|\beta|}{2}\right) \sim
\frac{\sinh\frac{L_0}{2}\sinh\frac{L_j}{2}}{\cosh\frac{|\beta|}{2}}
\sim {\rm const.}\exp\left(-\frac{|\beta|}{2}\right)
\end{eqnarray*}
as $|\beta| \rightarrow \infty$.
\end{pf}

\vskip 10pt

\noindent {\bf Geometric interpretation} \,\, We would like to
explore the geometric meanings of the summands in the complexified
formula (\ref{eqn:reform of non-cusp cases with GS}).

\vskip 10pt

In the case that $M$ has no cone points, all its geometric
boundary components (here cusps are treated as boundary geodesics
of length $0$) $\Delta_0, \Delta_1, \cdots, \Delta_n$ are boundary
geodesics with hyperbolic lengths $L_0, L_1, \cdots, L_n$
respectively. Assume $\Delta_0$ is not a cusp, that is, $L_0>0$.
Then as explained in \S \ref{s:calculation}, in the first sum the
summand is the width of one of the main gaps in the pair of pants
$\mathcal P(\Delta_0,\alpha,\beta)$ bounded by $\Delta_0$ and
$\alpha, \beta$; while in the second sum the sub-summand is the
width of one of the two extra gaps associated to $\Delta_j$ in the
pair of pants $\mathcal P(\Delta_0,\Delta_j,\beta)$ bounded by
$\Delta_0, \Delta_j$ and $\beta$. We would like to think of the
union of the two extra gaps in $\mathcal
P(\Delta_0,\Delta_j,\beta)$ as the orthogonal projection of
$\Delta_i$ onto $\Delta_0$ along the common perpendicular $\delta$
of $\Delta_j$ and $\Delta_0$ in $\mathcal
P(\Delta_0,\Delta_j,\beta)$ and think of its width as the {\it
direct visual measure} of $\Delta_j$ at $\Delta_0$ along $\delta$.
Hence the second part of the left hand side of (\ref{eqn:reform of
non-cusp cases with GS}) can be thought of as the total direct
visual measure of all the non-distinguished geometric boundary
components $\Delta_1, \cdots, \Delta_n$ at $\Delta_0$.

\vskip 10pt

In the case that $\Delta_0$ is a cone point of angle $\theta_0 \in
(0, \pi]$ (hence $L_0=\theta_0 i$) and all other geometric
boundary components of $M$ are boundary geodesics (here cusps
treated as boundary geodesics of length $0$), for each pair of
generalized simple closed geodesics $\alpha, \beta$ which bound
with $\Delta_0$ an embedded pair of pants $\mathcal
P(\Delta_0,\alpha,\beta)$ on $M$, each of $\alpha, \beta$ has a
direct visual angle at the cone point $\Delta_0$; and the summand
in the first sum is $i$ times the angle measure of one of the two
gaps at $\Delta_0$ between the two $\Delta_0$-geodesic rays
asymptotic to $\alpha^+, \beta^-$(respectively $\alpha^-$,
$\beta^+$ ). The sub-summand in the second sum is $i$ times half
the visual angle measure of $\Delta_j$ at $\Delta_0$ in the pair
of pants $\mathcal P(\Delta_0,\Delta_j,\beta)$ on $M$.

\vskip 10pt

When $M$ has cone points other than $\Delta_0$, the similar
formulations of the generalized McShane's identities
(\ref{eqn:001})--(\ref{eqn:00'}) in terms of ${\rm Gap}(\Delta_0;
\alpha, \beta)$ will not be as neat as in the above two special
cases. The problem lies in that a cone point (other than
$\Delta_0$) {\it seems} to have direct visual measure zero at
$\Delta_0$, causing the formulas to be non-uniform. However, this
non-uniformity is caused by the (wrong) point of view that we
treat a cone point as only a point. The correct point of view is
(perhaps) that a cone point (as a geometric boundary component)
should be a geodesic perpendicular to the surface at the very cone
point when the surface is ``imagined'' as lying in the hyperbolic
3-space and hence one should use purely complex length instead of
real one for a cone point. (The point of view of using complex
translation length for an isometry of the hyperbolic 3-space is
well discussed in details in \cite{fenchel1989book} and
\cite{series2001pjm}.)

\vskip 10pt

First assume that $\Delta_0$ is boundary geodesic of length $l_0
>0$ and consider a pair of generalized simple closed geodesics
$\alpha, \beta$ on $M$ such that $\alpha$ is a cone point of angle
$\varphi \in (0,\pi]$ and $\beta$ is an interior generalized
simple closed geodesic and that they bound with $\Delta_0$ an
embedded pair of pants $\mathcal P(\Delta_0,\alpha,\beta)$ on $M$.

Let the (unoriented) geodesic arc in $\mathcal
P(\Delta_0,\alpha,\beta)$ which is perpendicular to $\Delta_0$ and
$\alpha$ (respectively, $\alpha$ and $\beta$, $\beta$ and
$\Delta_0$) be denoted $[\Delta_0,\alpha]$ (respectively,
$[\alpha, \beta]$, $[\beta,\Delta_0]$). We cut $\mathcal
P(\Delta_0,\alpha,\beta)$ open along $[\Delta_0,\alpha]$,
$[\alpha, \beta]$, $[\beta,\Delta_0]$ to obtain two congruent
pentagons; lift one of them to a pentagon ${\mathbf
P}(\Delta_0,\alpha,\beta)$ in the hyperbolic plane $H^2$. Then by
Fenchel \cite{fenchel1989book} ${\mathbf
P}(\Delta_0,\alpha,\beta)$ can be regarded as a right angled
hexagon ${\mathbf H}(\Delta_0,\tilde{\alpha},\beta)$ spanned by
straight lines $\Delta_0,\tilde{\alpha},\beta$ in a hyperbolic
3-space $H^3$ containing the hyperbolic plane $H^2$. See Figure
\ref{fig08} for an illustration. Here $\tilde{\alpha}$ is the
straight line in $H^3$ which passes through the cone point
$\alpha$ in $H^2$ and is perpendicular to $H^2$. Let the common
perpendiculars in $H^3$ between pairs of
$\Delta_0,\tilde{\alpha},\beta$ be $[\Delta_0,\tilde{\alpha}]$,
$[\tilde{\alpha}, \beta]$, $[\beta,\Delta_0]$, where, as straight
lines $[\Delta_0,\tilde{\alpha}]$, $[\tilde{\alpha}, \beta]$ are
the same as $[\Delta_0,\alpha]$, $[\alpha, \beta]$ respectively.

We orient the six straight lines in the cyclic order
$\Delta_0,[\Delta_0,\tilde{\alpha}],\tilde{\alpha},[\tilde{\alpha},
\beta],\beta,[\beta,\Delta_0]$ as Fenchel did in
\cite{fenchel1989book}; see Figure \ref{fig08}. Then the three
oriented sides $\Delta_0,\tilde{\alpha},\beta$ of the right angled
hexagon ${\mathbf H}(\Delta_0,\tilde{\alpha},\beta)$ have complex
lengths $\frac{l_0}{2} + \pi i, \frac{\varphi i}{2} + \pi i,
\frac{|\beta|}{2} + \pi i$ respectively.

Let the ideal points which are the starting and ending endpoints
of an oriented straight line ${\bf l}$ in $H^3$ be denoted ${\bf
l}^{-},{\bf l}^{+}$ respectively. Then we have in $H^3$ an
oriented straight line $[\Delta_0,\tilde{\alpha}^{+}]$ which
intersects $\Delta_0$ perpendicularly and has $\tilde{\alpha}^{+}$
as its ending ideal point, and similarly an oriented straight line
$[\Delta_0, \beta^{-}]$ which intersects $\Delta_0$
perpendicularly and has $\beta^{-}$ as its ending ideal point.

\vskip 8pt

Then it can be verified that
$$
G\left(\frac{l_0}{2}, \frac{\varphi i}{2}+\pi i, \frac{|\beta|}{2}
+ \pi i\right)=G\left(\frac{l_0}{2}, \frac{\varphi i}{2},
\frac{|\beta|}{2}\right)
$$
is the complex length from $[\Delta_0, \beta^{-}]$ to
$[\Delta_0,\tilde{\alpha}^{+}]$ measured along $\Delta_0$ and
$$
S\left(\frac{l_0}{2}, \frac{\varphi i}{2}+\pi i, \frac{|\beta|}{2}
+ \pi i\right)=S\left(\frac{l_0}{2}, \frac{\varphi i}{2},
\frac{|\beta|}{2}\right)
$$
is the the complex length from $[\Delta_0,\tilde{\alpha}^{+}]$ to
$[\Delta_0,\tilde{\alpha}]$ measured along $\Delta_0$.

Note that $S(\frac{l_0}{2}, \frac{\varphi i}{2},
\frac{|\beta|}{2})$ is purely imaginary, which is obvious from its
geometric meaning.

\vskip 8pt

\begin{rmk}
We remark that it is crucial that in $G(\frac{l_0}{2},
\frac{\varphi i}{2}+\pi i, \frac{|\beta|}{2} + \pi i)$ and
\newline $S(\frac{l_0}{2}, \frac{\varphi i}{2}+\pi i, \frac{|\beta|}{2}
+ \pi i)$ the value used for $\Delta_0$ is $\frac{l_0}{2}$ instead
of $\frac{l_0}{2}+\pi i$.
\end{rmk}

\vskip 10pt

\begin{figure}
\begin{center}
\mbox{\beginpicture \setcoordinatesystem units <0.035in,0.035in>
\setplotarea x from 0 to 40, y from -25 to 22 \setdashes<2pt>
\arrow <6pt> [.16,.6] from 0 15 to 0 -17 \setquadratic \plot 24 20
9 1 2 -17 / \setsolid \arrow <6pt> [.16,.6] from 24.9 -14.92 to 25
-15 \arrow <6pt> [.16,.6] from 2.03 -16.8 to 2 -17 \arrow <6pt>
[.16,.6] from 0 -14 to 0 -17 \arrow <6pt> [.16,.6] from 24 18 to
-4 -3 \arrow <6pt> [.16,.6] from 40 5 to 15 15.75 \plot -3 1 10 -5
25 -15 / \plot 13 -20 25 -10 45 -3 / \arrow <6pt> [.16,.6] from
44.5 -3.1 to 45 -3 \plot 34 -10 35 2 43 15 / \arrow <6pt> [.16,.6]
from 42.92 14.9 to 43 15 \plot 38 17 20 0 11 -20 / \arrow <6pt>
[.16,.6] from 11.2 -19 to 11 -20 \put {\mbox{\scriptsize
$[\triangle_0,\alpha]$ }} [rt] <-1mm,0mm> at -4 -3 \put
{\mbox{\scriptsize $[\beta,\triangle_{0}]$ }} [ct] <1mm,0mm> at 34
-10 \put {\mbox{\scriptsize $[\alpha,\beta]$ }} [ct] <1mm,0mm> at
25 -15 \put {\mbox{\scriptsize $[\triangle_0,\beta^{-}]$ }} [cb]
<1mm,1mm> at 38 17 \put {\mbox{\scriptsize
$[\triangle_0,\tilde{\alpha}^{+}]$ }} [cb] <0mm,1mm> at 24 20 \put
{\mbox{\scriptsize $\triangle_{0}$ }} [rc] <0mm,0mm> at 15 17 \put
{\mbox{\scriptsize $\beta$ }} [lc] <0.7mm,0mm> at 45 -3 \put
{\mbox{\scriptsize $\tilde{\alpha}$ }} [cb] <1mm,1mm> at 0 15 \put
{\mbox{\scriptsize $\alpha$ }} [rb] <0mm,1mm> at 0 1
\ellipticalarc axes ratio 2:1.3 -80 degrees from 1.2 0.8 center at
0 0 \setlinear \plot 20 13.5 19 12.8 18 13.2 / \plot 30 9.4 28.8
8.4 28 8.8 / \plot 36 6.5 35.3 5.5 36 5.1 / \plot 33.6 -5.2 32
-5.8 32 -6.8 / \plot 20.2 -11.5 21 -11 22 -11.7 / \put
{\mbox{\small \sc Figure 8 \label{fig08}}} [ct] <0mm,-8mm> at 20
-18
\endpicture}\
\hspace{0.3in} \raisebox{2ex}{\mbox{\beginpicture
\setcoordinatesystem units <0.035in,0.035in> \setplotarea x from 0
to 40, y from -23 to 23 \setdashes<2pt> \arrow <6pt> [.16,.6] from
0 15 to 0 -22 \arrow <6pt> [.16,.6] from 18.7 22 to 18.7 3
\setquadratic \plot 24 15 9 -4 2 -22 / \setsolid \arrow <6pt>
[.16,.6] from 2.03 -21.8 to 2 -22 \arrow <6pt> [.16,.6] from 0 -20
to 0 -22 \arrow <6pt> [.16,.6] from 24 18 to -4 -3
\arrow <6pt> [.16,.6] from 40 5 to 10 17.9 \plot -3 0 10 -2 28 -9
/ \arrow <6pt> [.16,.6] from 27.8 -8.9 to 28 -9 \plot 12 -15 22 0
38 14 /
\plot 23 20.5 16 8 9.7 -15 / \arrow <6pt> [.16,.6] from 37.9 13.93
to 38 14 \arrow <6pt> [.16,.6] from 9.72 -14.9 to 9.7 -15 \arrow
<6pt> [.16,.6] from 18.7 5 to 18.7 3 \setlinear
\plot 18.7 10.3 18 9.5  18 8.5 / \plot 17.2 -4.5 18 -3.55 18.7
-3.8 / \plot 30.3 9.1 29.4 8.3 30.2 8. / \circulararc -40 degrees
from 1.2 0.8  center at 0 0 \put {\mbox{\LARGE $\cdot$}} [cc]
<0mm,0mm> at 18.9 13.9 \put {\mbox{\LARGE $\cdot$}} [cc] <0mm,0mm>
at 0. -0.2 \ellipticalarc axes ratio 2:1.3 120 degrees from 18
13.3 center at 19.5 13.8 \put {\mbox{\scriptsize
$[\triangle_0,\alpha]$}} [rt] <-1mm,0mm> at 0.2 -4 \put
{\mbox{\scriptsize $[\beta,\triangle_{0}]$}} [lt] <-1mm,-1mm> at
40 5 \put {\mbox{\scriptsize $[\alpha,\beta]$}} [ct] <1mm,0mm> at
27.8 -8.9 \put {\mbox{\scriptsize
$[\tilde{\triangle}_0,\beta^{-}]$ }} [lb] <1mm,1mm> at 23 20.5
\put {\mbox{\scriptsize $[\triangle_0,\tilde{\alpha}^{+}]$ }} [ct]
<1mm,-1mm> at 7 -22 \put {\mbox{\scriptsize $\triangle_{0}$ }}
[rc] <-3mm,0mm> at 20 14 \put {\mbox{\scriptsize $\beta$}} [lc]
<1mm,1mm> at 38 14 \put {\mbox{\scriptsize $\tilde{\triangle}_{0}$
}} [cb] <1mm,1mm> at 18.7 22 \put {\mbox{\scriptsize
$\tilde{\alpha}$ }} [cb] <1mm,1mm> at 0 15 \put {\mbox{\scriptsize
$\alpha$}} [rb] <-1mm,1mm> at 0.7 -0.2 \put {\mbox{\small \sc
Figure 9 \label{fig09}}} [ct] <0mm,-8mm> at 17 -22
\endpicture}}
\end{center}
\end{figure}

Next assume $\Delta_0$ is a cone point of angle $\theta_0 \in (0,
\pi]$ and consider a pair of generalized simple closed geodesics
$\alpha, \beta$ on $M$ such that $\alpha$ is a cone point of angle
$\varphi \in (0,\pi]$ and $\beta$ is an interior generalized
simple closed geodesic and that they bound with $\Delta_0$ an
embedded pair of pants $\mathcal P(\Delta_0,\alpha,\beta)$ on $M$.

In this case we cut $\mathcal P(\Delta_0,\alpha,\beta)$ open along
$[\Delta_0,\alpha]$, $[\alpha, \beta]$, $[\beta,\Delta_0]$ to
obtain two congruent quadrilaterals and lift one of them to a
quadrilateral ${\mathbf Q}(\Delta_0,\alpha,\beta)$ in the
hyperbolic plane $H^2$. As before, let $\tilde{\alpha}$ be the
straight line in $H^3$ which passes the cone point $\alpha$ in
$H^2$ and is perpendicular to $H^2$. Similarly for
$\tilde{\Delta_0}$. Then we obtain a right angled hexagon
${\mathbf H}(\tilde{\Delta_0},\tilde{\alpha},\beta)$ in $H^3$. Let
the six sides of ${\mathbf
H}(\tilde{\Delta_0},\tilde{\alpha},\beta)$ be oriented as
illustrated in Figure \ref{fig09}. Then the three oriented sides
$\tilde{\Delta_0},\tilde{\alpha},\beta$ of the right angled
hexagon ${\mathbf H}(\tilde{\Delta_0},\tilde{\alpha},\beta)$ have
complex lengths $\frac{\theta_0 i}{2} + \pi i, \frac{\varphi}{2}i
+ \pi i, \frac{|\beta|}{2} + \pi i$ respectively.

Similarly, we have in $H^3$ an oriented straight line
$[\tilde{\Delta_0},\tilde{\alpha}^{+}]$ which intersects
$\tilde{\Delta_0}$ perpendicularly and has $\tilde{\alpha}^{+}$ as
its ending ideal point, and another oriented straight line
$[\tilde{\Delta_0}, \beta^{-}]$ which intersects
$\tilde{\Delta_0}$ perpendicularly and has $\beta^{-}$ as its
ending ideal point.

Then
$$
G\left(\frac{\theta_0 i}{2}, \frac{\varphi i}{2}+\pi i,
\frac{|\beta|}{2} + \pi i\right)=G\left(\frac{\theta_0 i}{2},
\frac{\varphi i}{2}, \frac{|\beta|}{2}\right)
$$
is the complex length from $[\tilde{\Delta_0}, \beta^{-}]$ to
$[\tilde{\Delta_0},\tilde{\alpha}^{+}]$ measured along
$\tilde{\Delta_0}$ and
$$
S\left(\frac{\theta_0 i}{2}, \frac{\varphi i}{2}+\pi i,
\frac{|\beta|}{2} + \pi i\right)=S\left(\frac{\theta_0 i}{2},
\frac{\varphi i}{2}, \frac{|\beta|}{2}\right)
$$
is the the complex length from
$[\tilde{\Delta_0},\tilde{\alpha}^{+}]$ to
$[\tilde{\Delta_0},\tilde{\alpha}]$ measured along
$\tilde{\Delta_0}$.

Note that $S(\frac{\theta_0 i}{2}, \frac{\varphi i}{2}+\pi i,
\frac{|\beta|}{2} + \pi i)$ is real, which is obvious from its
geometric meaning.

\vskip 8pt

\begin{rmk}
Here it is crucial that in $G(\frac{\theta_0 i}{2}, \frac{\varphi
i}{2}+\pi i, \frac{|\beta|}{2} + \pi i)$ and
\newline $S(\frac{\theta_0 i}{2}, \frac{\varphi i}{2}+\pi i, \frac{|\beta|}{2}
+ \pi i)$ the value used for $\Delta_0$ is $\frac{\theta_0 i}{2}$
instead of $\frac{\theta_0 i}{2}+\pi i$.
\end{rmk}

\end{document}